\documentclass[12pt]{article}
\usepackage{fullpage,graphicx,psfrag,url,subfigure}
\usepackage[small,bf]{caption}
\usepackage{pstricks}
\newcommand{\devices}{\mathcal{D}}
\newcommand{\terminals}{\mathcal{T}}
\newcommand{\nets}{\mathcal{N}}

\newcommand{\eg}{{\it e.g.}}
\newcommand{\ie}{{\it i.e.}}

\newcommand{\reals}{{\mbox{\bf R}}}
\newcommand{\argmin}{\mathop{\rm argmin}}

\title{Message Passing for Dynamic Network Energy Management} 
\author{Matt Kraning, Eric Chu, Javad Lavaei, and
Stephen Boyd\thanks{All authors are members of the Information Systems
Laboratory, Department of Electrical Engineering, Stanford University}
\date{WORKING DRAFT --- \today}}

\begin{document}
\maketitle

\begin{abstract}
We consider a network of devices, such as generators, fixed loads,
deferrable loads, and storage devices, each with its own dynamic
constraints and objective, connected by lossy capacitated lines. The
problem is to minimize the total network objective subject to the device
and line constraints, over a given time horizon. This is a large
optimization problem, with variables for consumption or generation
in each time period for each device.  In this paper we
develop a decentralized method for solving this problem. The method is
iterative: At each step, each device exchanges simple messages with its
neighbors in the network and then solves its own optimization problem,
minimizing its own objective function, augmented by a term determined by
the messages it has received. We show that this message passing method
converges to a solution when the device objective and constraints are
convex. The method is completely decentralized, and needs no global
coordination other than synchronizing iterations; the problems to be
solved by each device can typically be solved extremely efficiently and
in parallel. The method is fast enough that even a serial implementation
can solve substantial problems in reasonable time frames. We report
results for several numerical experiments, demonstrating the method's
speed and scaling, including the solution of a problem instance with
over $30$ million variables in $52$ minutes for a serial
implementation; with decentralized computing, the solve time would be
less than one second.
\end{abstract}

\newpage
\tableofcontents
\newpage

\section{Introduction} \label{s-intro}

A traditional power grid is operated by solving a number of optimization
problems. At the transmission level, these problems include unit
commitment, economic dispatch, optimal power flow (OPF), and
security-constrained OPF (SCOPF). At the distribution level, these
problems include loss minimization and reactive power compensation. With
the exception of the SCOPF, these optimization problems are static with
a modest number of variables (often less than $10000$), and are solved
on time scales of $5$ minutes or more. However, the operation of next
generation electric grids (\ie, smart grid) will rely critically on solving
large-scale, dynamic optimization problems involving hundreds of thousands of
devices jointly optimizing tens to hundreds of
millions of variables, on the order of seconds rather than minutes
\cite{ET:11,LVH:12}. More precisely, the distribution level of a smart
grid will include various types of active dynamic devices, such as
distributed generators based on solar and wind, batteries, deferrable
loads, curtailable loads, and electric vehicles, whose control and
scheduling amount to a very complex power management problem
\cite{TSBC:11,CI:12}.

In this paper, we consider a general problem, which we call the
\emph{optimal power scheduling problem} (OPSP), in which a network of
dynamic devices are connected by lossy capacitated lines, and the goal is to
jointly minimize a network objective subject to local constraints on the
devices and lines. The network objective is the sum of the objective
functions of each device. These objective functions extend over a given time
horizon and encode operating costs such as fuel consumption and constraints
such as limits on power generation or consumption.  In addition, the objective
functions encode dynamic objectives and constraints such as limits on ramp
rates or charging limits. The variables for each device consist of its
consumption or generation in each time period and can also include local
variables which represent internal states of the device over time, such as the
state of charge of a storage device.

When all device objective functions and line constraints are convex, the
OPSP is a convex optimization problem, which can in principle be solved
efficiently \cite{BoV:04}. If not all device objective functions are
convex, we can solve a relaxed form of the OPSP which can be used to
find good, local solutions to the OPSP. The optimal value of the relaxed
OPSP also gives a lower bound for the optimal value of the OPSP which
can be used to evaluate the suboptimality of a local solution.

For any network, the corresponding OPSP contains at least as many
variables as the number of devices multiplied by the length of the time
horizon. For large networks with hundreds of thousands of devices and a
time horizon with tens or hundreds of time periods, the extremely large
number of variables present in the corresponding OPSP makes solving it
in a centralized fashion computationally impractical, even when all
device objective functions are convex.

We propose a decentralized optimization method which efficiently solves
the OPSP by distributing computation across every device in the network.
This method, which we call \emph{prox-average message passing}, is
iterative: At each iteration, every device passes simple messages to its
network neighbors and then solves an optimization problem that minimizes
the sum of its own objective function and a simple regularization term
that only depends on the messages it received from its network neighbors
in the previous iteration. As a result, the only non-local coordination
needed between devices for prox-average message passing is 
synchronizing iterations. When all device objective functions are
convex, we show that prox-average message passing converges to a
solution of the OPSP.

Our algorithm can be used several ways.  
It can be implemented in a traditional way on a 
single computer or cluster of computers, 
by collecting all the device constraints and objectives.
We will demonstrate this use with an implementation that runs on a single
8-core computer.
A more interesting use is in a peer-to-peer architecture.
In this architecture, each device contains its own processor, which carries 
out the required local dynamic optimization and exchanges messages with its
neighbors on the network.
In this setting, the devices do not need to divulge their objectives or 
constraints; they must only support a simple protocol for interacting
with the neighbors.
Our algorithm ensures that the network power flows will converge to their
optimal values, even though each device has very little information about
the rest of the network, and only exchanges limited messages with its
immediate neighbors.

Due to recent advances in convex optimization \cite{WB:10, CVXGEN,
MWB:11}, in many cases the optimization problems that each device solves
in each iteration of prox-average message passing can be executed at
millisecond or even microsecond time-scales on inexpensive, embedded
processors. Since this execution can happen in parallel across all
devices, the entire network can execute prox-average message passing at
kilohertz rates. We present a series of numerical examples to illustrate
this fact by using prox-average message passing to solve instances of
the OPSP with up to $10$ million variables serially in $17$ minutes.
(This is on an 8-core computer; with 64 cores, the time would be around
$2$ minutes.)
Using decentralized computing, the solve time would be essentially
independent of the size of the network and measured in fractions of a second.

We note that although the primary application for our method is power
management, it can easily be adapted to more general resource allocation
and graph-structured optimization problems \cite{RAW:10}.

\paragraph{Related work.}

The use of optimization in power systems dates back to the 1920s and has
traditionally concerned the optimal dispatch problem \cite{Happ:77},
which aims to find the lowest cost method for generating and delivering
power to consumers, subject to physical generator constraints. With the
advent of computer and communication networks, many different ways to
numerically solve this problem have been proposed \cite{Zhu:09} and more
sophisticated variants of optimal dispatch have been introduced, such as
OPF, economic dispatch, and dynamic dispatch \cite{CR:90}, which extend
optimal dispatch to include various reliability and dynamic constraints.
For reviews of optimal and economic dispatch as well as general power
systems, we direct the reader to \cite{BV:99} and the book and review
papers cited above.

When modeling AC power flow, the OPSP is a dynamic version of the OPF
\cite{Carp:62}, extending the latter to include many more types of
devices such as storage units. The OPSP also introduces a time horizon with
coupling constraints between variables across time periods. The OPF has
been a fundamental problem in power systems for over $50$ years and is
known to be non-convex. Recently, it was shown that the OPF can be
solved exactly in certain circumstances by recasting it as a
semidefinite program and solving its dual problem \cite{LL:12}.

Distributed optimization methods are naturally applied to power networks
given the graph-structured nature of the transmission and
distribution networks. There is an extensive literature on distributed
optimization methods, dating back to the early 1960’s. The prototypical
example is dual decomposition \cite{DaW:60, Eve:63}, which is based on
solving the dual problem by a gradient method.
In each iteration, all devices optimize their local (primal)
variables based on current prices (dual variables). Then the dual
variables are updated to account for imbalances in supply and demand,
with the goal being to determine prices for which supply equals demand.

Unfortunately, dual decomposition methods are not robust, requiring many
technical conditions, such as strict convexity and finiteness of all
local cost functions, for both theoretical and practical convergence.
One way to loosen these technical conditions is to use an augmented
Lagrangian \cite{Hes:69a,Pow:69,Ber:82}, resulting in the method of
multipliers. This subtle change allows the method of multipliers to
converge under mild technical conditions, even when the local cost
functions are not strictly convex or necessarily finite. However, this
method has the disadvantage of no longer being separable across
subsystems. To achieve both separability and robustness for distributed
optimization, we can instead use the \emph{alternating direction method
of multipliers} (ADMM) \cite{GlM:75, GaM:76, Eck:94b, BP:11}.
ADMM is very closely related to many other algorithms, and is 
identical to Douglas-Rachford operator splitting; 
see, \eg, the discussion in \cite[\S 3.5]{BP:11}.

Building on the work of \cite{LL:12}, a distributed algorithm was
recently proposed \cite{LZT:11} to solve the dual OPF using a standard
dual decomposition on subsystems that are maximal cliques of the power
network. Augmented Lagrangian methods (including ADMM) have previously
been applied to the study of power systems with static,
single period objective functions on a small number of distributed
subsystems, each representing regional power generation and consumption
\cite{KB:97}. For an overview of related decomposition methods applied
to power flow problems, we direct the reader to \cite{KB:00, Bal:06} and
the references therein.

% Distributed optimization techniques date back to the early 1960's, with the
% most famous example being the dual decomposition \cite{Eve:63}.  In this
% framework, rather than solving the primal optimization problem directly, one
% forms the Lagrangian and iteratively updates the primal and dual variables to
% find a saddle point of the Lagrangian.  These updates can be conducted in
% parallel by having each device optimize their local variables based on the
% current dual variables, which can be interpreted as prices.  The dual
% variables are then updated in response to the new supply and demand given by
% the devices' new consumption and generation profiles.
% 
% The main drawbacks of dual decomposition are the large number of technical
% conditions needed to guarantee both theoretical and practical convergence, in
% particular the requirements of strict convexity and finiteness for all device
% objective functions.  Augmented Lagrangian methods
% \cite{Hes:69a,Pow:69,Ber:82} were developed as a way to bring more robustness
% to the dual decomposition and loosen the technical conditions required for
% convergence.  By replacing the Lagrangian with the augmented Lagrangian, the
% resulting algorithm, called the \emhp{method of multipliers}, achieves robust
% convergence guarantees without requiring strict convexity or finiteness of
% device objective function, but at the expense of destroying the separability
% of the optimization across devices.

\paragraph{Outline.} The rest of this paper is organized as follows. In \S
\ref{s-model} we give the formal definition of our network model. In
\S \ref{s-devices} we give examples of how to model specific
devices such as generators, deferrable loads and energy storage systems
in our formal framework. In \S \ref{s-convexity}, we describe the
role that convexity plays in the OPSP and introduce the idea of convex
relaxations as a tool to find solutions to the OPSP in the presence of
non-convex device objective functions. In \S \ref{s-method} we
derive the prox-average message passing equations. In \S
\ref{s-numerical} we present a series of numerical examples, and in
\S \ref{s-extensions} we discuss how our framework can be easily
extended to include use cases we do not explicitly cover in this paper.

\section{Network model} \label{s-model}

\subsection{Formal definition}
A network consists of a finite set of \emph{terminals} $\terminals$, a
finite set of \emph{devices} $\devices$, and a finite set of \emph{nets}
$\nets$. The sets $\devices$ and $\nets$ are both partitions of
$\terminals$. Thus, each terminal is associated with exactly one device
and exactly one net. Equivalently, a network can be defined as a
bipartite graph with one set of vertices given by devices, the other set
of vertices given by nets, and edges given by terminals.
%We can
%equivalently consider a network as a bipartite graph, with one set of
%nodes given by the devices, the other set of nodes given by the nets,
%and terminals acting as edges between those two sets.

Each terminal $t\in\terminals$ has an associated \emph{power schedule}
$p_t = (p_t(1), \ldots, p_t(T)) \in \reals^T$, where $T$ is a given time
horizon.  Here $p_t(\tau)$ is the amount of energy
consumed by device $d$ through terminal $t$ in time period $\tau$, where
$t\in d$ (\ie, terminal $t$ is associated with device $d$).  When
$p_t(\tau) < 0$, $-p_t(\tau)$ is the energy generated by device $d$
through terminal $t$ in time period $\tau$.  The set of all terminal
power schedules is denoted $p$.  This is a function from $\terminals$
(the set of terminals) into $\reals^T$ (time periods); we can identify
$p$ with a $|\terminals| \times T$ matrix, whose rows are the terminal
power schedules.

%For any subset of terminals $r \subseteq \terminals$, we use the notation
%$p_r$ to denote a vector with components $p_t$, for $t\in r$.

For each device $d\in\devices$, $p_d$ consists of the set of power
schedules of terminals associated with $d$, which we identify with a
$|d| \times T$ matrix whose rows are taken from the rows of $p$
corresponding to the terminals in $d$. Each device $d$ has an associated
\emph{objective function} $f_d:\reals^{|d|\times T} \to \reals \cup
\{+\infty\}$, where we set $f_d(p_d) = \infty$ to encode constraints on
the power schedules for the device. When $f_d(p_d)<\infty$, we say that
$p_d$ is a set of \emph{realizable} power schedules for device $d$, and
we interpret $f_d(p_d)$ as the cost (or revenue, if negative) to device
$d$ for the power schedules $p_d$.
 
Similarly, for each net $n \in \nets$, $p_n$ consists of the set of
power schedules of terminals associated with $n$, which we identify with
a $|n| \times T$ matrix whose rows are taken from the rows of $p$
corresponding to the terminals in $n$. Each net $n$ is a lossless energy
carrier (commonly referred to as a bus in power systems literature),
which is represented by the constraint \begin{equation}
\label{e-net-balance} \sum_{t \in n} p_t(\tau)=0, \quad \tau=1, \ldots,
T. \end{equation} In other words, in each time period the power flows on
each net balance.

For any terminal, we define the \emph{average net power imbalance} $\bar
p:\terminals \to \reals^T$, as
\begin{equation} \label{e-bar}
\bar p_t = \frac{1}{|n|}\sum_{t' \in n} p_{t'},
\end{equation}
where $t\in n$, \ie, terminal $t$ is associated with net $n$. In other
words, $\bar p_t(\tau)$ is the average power schedule of all terminals
associated with the same net as terminal $t$ at time $\tau$. We overload
this notation for devices
%%%%%%%%%as the average power schedules for all terminals on the device
by defining $\bar{p}_d = (\bar{p}_{i_1}, \ldots, \bar{p}_{i_{|d|}})$,
where device $d$ is associated with terminals $i_1, \ldots, i_{|d|}$.
Using an identical notation for nets, we can see that $\bar{p}_n$ simply
contains $|n|$ copies of the average net power imbalance for net $n$.
The net power balance constraint for all terminals can be equivalently
expressed as $\bar p =0$.

\paragraph{Optimal power scheduling problem.}
We say that a set of power schedules $p: \terminals \to \reals^T$ is
\emph{feasible} if $f_d(p_d) < \infty $ for all $d \in \devices$ (\ie,
all devices' power schedules are realizable), and $\bar p = 0$ (\ie,
power balance holds across all nets). We define the network objective as
$f(p) = \sum_{d \in \devices}f_d(p_d)$. The \emph{optimal power
scheduling problem} (OPSP) is
\begin{equation}\label{e-network}
\begin{array}{ll}
\mbox{minimize} & f(p)\\
\mbox{subject to} & \bar{p} = 0,
\end{array}
\end{equation}
with variable $p: \terminals \to \reals^T$. We refer to $p$ as
\emph{optimal} if it solves~(\ref{e-network}), \ie, globally minimizes
the objective among all feasible $p$. We refer to $p$ as \emph{locally
optimal} if it is a locally optimal point for~(\ref{e-network}).

\paragraph{Dual variables and locational marginal prices.}
Suppose $p^0$ is a set of optimal power schedules, that also
minimizes the Lagrangian
\[
f(p) + \sum_{t\in\terminals} \sum_{\tau=1}^T y_t^0(\tau) \bar{p}_t(\tau),
\]
with no power balance constraint, where $y^0: \terminals \to \reals^T$.
In this case we call $y^0$ a set of optimal
Lagrange multipliers or dual variables.
When $p^0$ is a locally optimal point, which also locally minimizes
the Lagrangian, then we refer to $y^0$ as a set of locally optimal Lagrange
multipliers.

The dual variables $y^0$ are
related to the traditional concept of locational marginal prices
$\mathcal{L}^0:\terminals \to \reals^T$ by rescaling the dual variables
associated with each terminal according to the size of its associated
net, \ie, $\mathcal{L}_t^0 = y_t^0/|n|$, where $t\in n$. This rescaling
is due to the fact that locational marginal prices are the dual variables
associated with the constraints in (\ref{e-net-balance}) rather than
their scaled form used in (\ref{e-network}) \cite{EW:02}.

\subsection{Discussion}
We now describe our model in a more intuitive, less formal manner.
Devices include generators, loads, energy storage systems, and other
power sources, sinks, and converters. Terminals are ports on a device
through which power flows, either into or out of the device (or both,
at different times, as in a storage device).
Nets are used to model ideal lossless uncapacitated connections 
between terminals over which power
is transmitted; losses, capacities, and more general connection
constraints between a set of terminals can be modeled with the addition
of a device and individual nets which connect each terminal to the new
device. Our network model does not specify nor require a specific type
of energy transport mechanism (\eg, DC, single or 3-phase AC), but rather can
abstractly model arbitrary heterogeneous energy transport and exchange
mechanisms.

The objective function of a device is used to measure the cost (which
can be negative, representing revenue) associated with a particular mode
of operation, such as a given level of consumption or generation of 
power.  This cost can include the actual direct cost of operating
according to the given power schedules, such as a fuel cost,
as well as other costs such as
$\mathrm{CO}_2$ generation, or costs associated with increased maintenance
and decreased system lifetime due to 
structural fatigue. The objective function
can also include local variables other than power schedules, such as the state
of charge of a storage device.  

Constraints on the power schedules and internal variables for a device
are encoded by setting the objective function to $+\infty$ for power 
schedules that violate the constraints.
In many cases, a device's
objective function will only take on the values $0$ and $+\infty$, indicating
no local preference among feasible power schedules.

\subsection{Example transformation to abstract network model}
\begin{figure}
\begin{center}
\subfigure{
    \includegraphics[width=0.5\textwidth]{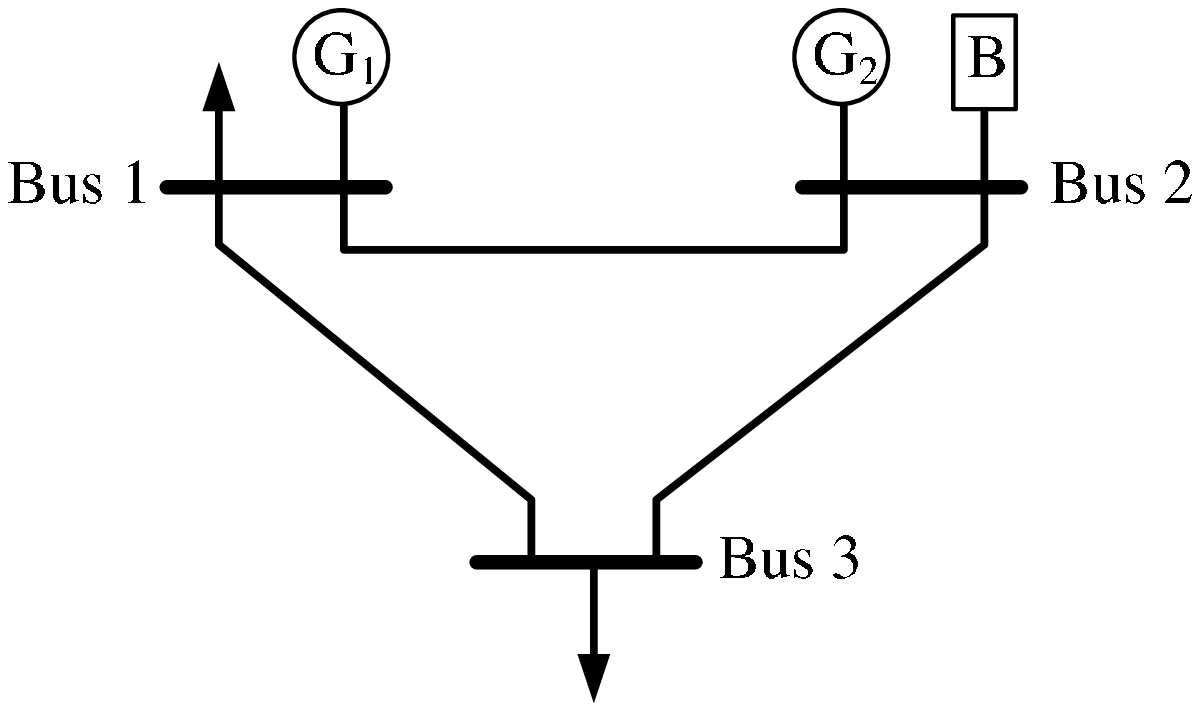}
            }
    \hspace{0.75cm}
\subfigure{
        \includegraphics[width=0.4\textwidth]{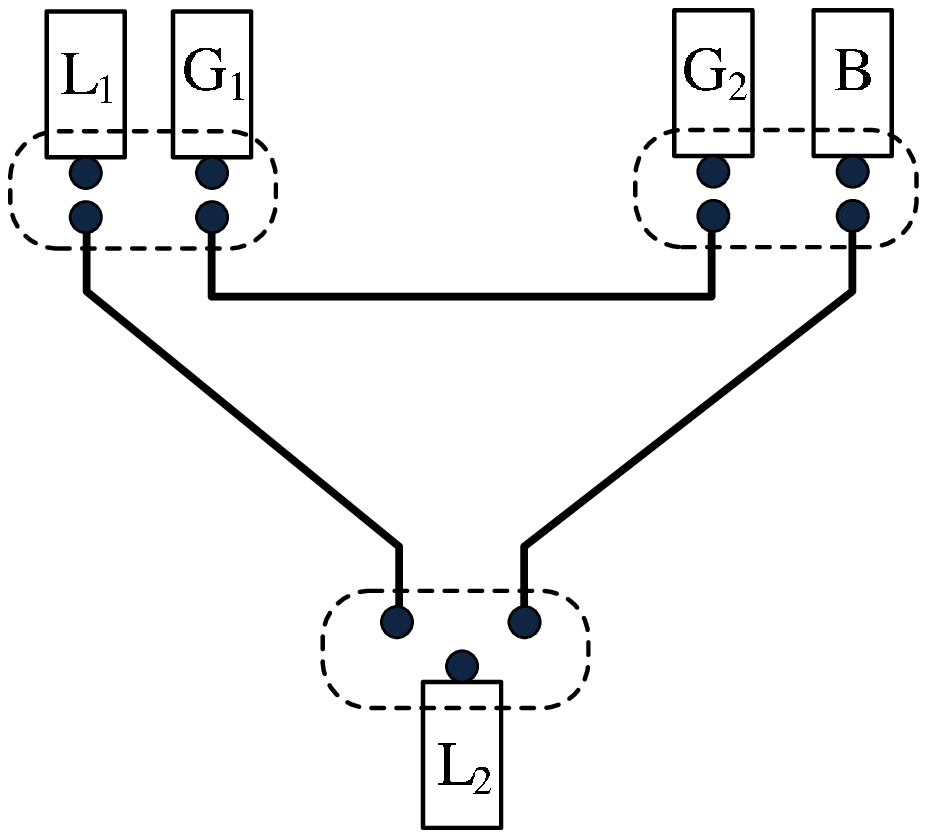}
            }
\end{center}
\caption{A simple network (\emph{left}),
and its transformation into standard form (\emph{right}).}

\label{f-network}
\end{figure}

We illustrate how a traditional power network can be recast into our
network model in figure \ref{f-network}. The original power network,
shown on the left, contains $2$ loads, $3$ buses, $2$ generators, and
a single battery storage system. We can
transform this small power grid into our model by representing it as a
network with $11$ terminals, $8$ devices ($3$ of them transmission
lines), and $3$ nets, shown on the right of figure~\ref{f-network}.
Terminals are shown as small filled circles.
Single terminal devices, which are used to model loads, generators,
and the battery, are shown as boxes.  The transmission lines are two terminal 
devices represented by solid lines.
The nets are shown as dashed rounded boxes. Terminals are associated with the
device they touch and the net in which they are contained.

\begin{figure}
\begin{center}
    \psfrag{1}[][]{$\mathrm{T}_3$}
    \psfrag{2}[][]{$\mathrm{L}_2$}
    \psfrag{3}[][]{$\mathrm{T}_2$}
    \psfrag{4}[][]{$\mathrm{G}_2$}
    \psfrag{5}[][]{$\mathrm{B}$}
    \psfrag{6}[][]{$\mathrm{T}_1$}
    \psfrag{7}[][]{$\mathrm{G}_1$}
    \psfrag{8}[][]{$\mathrm{L}_1$}
	\includegraphics[width=0.3\textwidth]{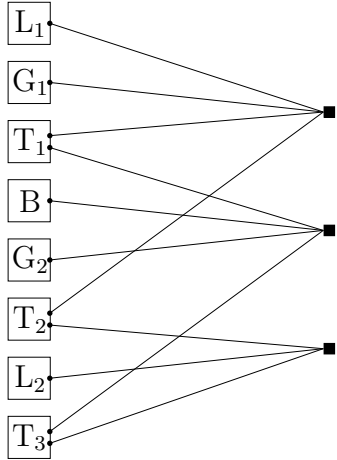}
\end{center}
\caption{The network in figure \ref{f-network} represented as a bipartite
graph. Devices (boxes) are shown on the left with their associated terminals (dots). 
The terminals are connected to their corresponding nets (solid boxes) on the right.}
\label{f-bipartite-graph}
\end{figure}

The set of terminals can be partitioned by either the devices 
they are associated with,
or the nets in which they are contained.
Figure \ref{f-bipartite-graph} shows the network in 
figure \ref{f-network} as a bipartite graph, with devices on the left and 
nets on the right. 

\section{Device examples}
\label{s-devices}
We present a few examples of how common devices can be modeled in
our framework. We note that these examples are intentionally kept simple,
but could easily be extended to more refined objectives and constraints. 
In these examples, it is easier to discuss operational
costs and constraints for each device separately. A device's objective
function is equal to the device's cost function unless any constraint is
violated, in which case we set the objective value to $+\infty$.

\paragraph{Generator.} A generator is a single-terminal device with
power schedule $p_\mathrm{gen}$, which generates power over a range,
$P^\mathrm{min} \leq -p_\mathrm{gen} \leq P^\mathrm{max}$, and has
ramp-rate constraints
\[
R^\mathrm{min} \leq -Dp_\mathrm{gen} \leq R^\mathrm{max}, 
\]
which limit the change of power levels from one period to the next.
Here, the operator $D\in\reals^{(T-1)\times T}$ is the forward
difference operator, defined as
\[
(Dx)(\tau) = x(\tau+1)-x(\tau), \quad \tau = 1,\ldots,T-1.
\]
The cost function for a generator has the separable form 
\[
\psi_\mathrm{gen}(p_\mathrm{gen}) =
\sum_{\tau=1}^T  \phi_\mathrm{gen}(-p_\mathrm{gen}(\tau)),
\]
where $\phi: \reals\to \reals$ gives the cost
of operating the generator at a given power level over a single time period.
This function is typically, but not always, convex and increasing. It could be
piecewise linear, or, for example, quadratic:
\[
\phi_\mathrm{gen}(u) = \alpha u^2 + \beta u,
\]
where $\alpha, \beta > 0$.  

More sophisticated models of generators allow for them to be switched on
or off, with an associated cost each time they are turned on or off.
When switched on, the generator operates as described above. When the
generator is turned off, it generates no power but can still incur costs
for other activities such as idling.

\paragraph{Transmission line.} A transmission line is a two terminal
device that transmits power across some distance with power schedules
$p_1$ and $p_2$ and zero cost function. The sum $p_1+p_2$ represents the
loss in the line and is always nonnegative. The difference $p_1-p_2$ can
be interpreted as twice the power flow from terminal one to terminal
two. A line has a maximum flow capacity, which is given by
\[
|p_1 - p_2| \leq C^\mathrm{max},
\]
as well as a loss function, $\ell(p_1,p_2):\reals^{2\times T} \to
\reals^T_+$, which defines the constraint
\[
p_1 + p_2 + \ell(p_1,p_2) = 0.
\]

In many cases, $\ell$ is a convex function with $\ell(0,0)=0$. Under a
simple resistive model, $\ell$ is a convex quadratic function of $p_1$
and $p_2$. Under a model for AC power transmission, the feasible region
defined by power loss is given by an ellipse \cite{LTZ:12}.

\paragraph{Battery.} A battery is a single terminal 
energy storage device with power schedule $p_\mathrm{bat}$, which can take in or
deliver energy, depending on whether it is charging or discharging. The
charging and discharging rates are limited by the constraints
$-D^\mathrm{max} \leq p_\mathrm{bat} \leq C^\mathrm{max}$, where
$C^\mathrm{max}\in\reals^T$ and $D^\mathrm{max}\in\reals^T$ are the
maximum charging and discharging rates. At time $t$, the charge level of
the battery is given by local variables
\[
q(\tau) = q^\mathrm{init}+\sum_{t=1}^\tau p_\mathrm{bat}(t),\quad \tau = 1,\ldots,T,
\] 
where $q^\mathrm{init}$ is the initial charge. It has zero cost function
and the charge level must not exceed the battery capacity, \ie, $0 \leq
q(\tau) \leq Q^\mathrm{max}$, $\tau = 1,\ldots,T$.
It is common to constrain the terminal battery charge $q(T)$ to be some 
specified value or to match the initial charge $q^\mathrm{init}$.

More sophisticated battery models include (possibly state-dependent)
charging and discharging inefficiencies as well as charge leakage. In
addition, they can include costs which penalize excessive charge-discharge
cycling.

\paragraph{Fixed load.} 
A fixed energy load is a single terminal device with zero cost function
which consists of a desired consumption profile, $l\in\reals^T$. This
consumption profile must be satisfied in each period, \ie, we have the
constraint $p_\mathrm{load} = l$.

\paragraph{Thermal load.} A thermal load is a single terminal device
with power schedule $p_\mathrm{therm}$ which consists of a heat store
(room, cooled water reservoir, refrigerator), with temperature profile
$\theta\in\reals^T$, which must be kept within minimum and maximum
temperature limits, $\theta^\mathrm{min}\in\reals^T$ and
$\theta^\mathrm{max}\in\reals^T$. The temperature of the heat store
evolves according to
\[
\theta(\tau+1) = \theta(\tau)+
	(\mu/c)(\theta^\mathrm{amb}(\tau)-\theta(\tau))-
	(\eta/c) p_\mathrm{therm}(\tau),
	\quad \tau = 1,\ldots,T-1,\qquad \theta(1) = \theta^\mathrm{init},
\]
where $0 \leq p_\mathrm{therm} \leq H^\mathrm{max}$ is the cooling power
consumption profile, $H^\mathrm{max}\in\reals^T$ is the maximum cooling
power, $\mu$ is the ambient conduction coefficient, $\eta$ is the
heating/cooling efficiency, $c$ is the heat capacity of the heat store,
$\theta^\mathrm{amb}\in\reals^T$ is the ambient temperature profile, and
$\theta^\mathrm{init}$ is the initial temperature of the heat store.
A thermal load has zero cost function.

More sophisticated models include temperature-dependent cooling and heating
efficiencies for heat pumps, more complex dynamics of the system
whose temperature is being controlled, and additional additive terms in
the thermal dynamics, to represent occupancy or other heat sources.

\paragraph{Deferrable load.} A deferrable load is a single terminal
device with zero cost function that must consume a minimum amount of
power over a given interval of time, which is characterized by the
constraint $\sum_{\tau=A}^{D} p_\mathrm{load}(\tau) \geq E$, where $E$ is the
minimum total consumption for the time interval $\tau = A,\ldots,D$. The energy
consumption in each time period is constrained by $0 \leq p_\mathrm{load} \leq
L^\mathrm{max}$. In some cases, the load can only be turned on or off in each
time period, \ie, $p_\mathrm{load}(\tau) \in \{0, L^\mathrm{max}\}$ for $\tau =
A, \ldots, D$.

\paragraph{Curtailable load.} A curtailable load is a single terminal
device which does not impose hard constraints on its power requirements,
but instead penalizes the shortfall between a desired load profile $l
\in \reals^T$ and delivered power. In the case of a linear penalty, its
cost function is given by
\[
\alpha (l - p_\mathrm{load})_+,
\]
where $(z)_+ = \max(0,z)$, $p_\mathrm{load} \in \reals^T$ is the amount
of electricity delivered to the device, and $\alpha > 0$ is a penalty
parameter.

\paragraph{Electric vehicle.} An electric vehicle is a single terminal
device with power schedule $p_\mathrm{ev}$ which has a desired charging
profile $c^\mathrm{des} \in \reals^T$ and can be charged within a time
interval $\tau = A,\ldots,D$. To avoid excessive charge cycling, we
assume that the electric vehicle battery cannot be discharged back into
the grid (in more sophisticated vehicle-to-grid models, this assumption
is relaxed), so we have the constraints $0 \leq p_\mathrm{ev} \leq
C^\mathrm{max}$, where $C^\mathrm{max}\in\reals^T$ is the maximum
charging rate. We assume that $c^\mathrm{des}(\tau) = 0$ for $\tau\notin
\{A,\ldots,D\}$, and for $\tau = A,\ldots,D$, the charge level is
given by
\[
q(\tau) = q^\mathrm{init}+\sum_{t = A}^\tau p_\mathrm{ev}(t),
\]
where $q^\mathrm{init}$ is the initial charge of the vehicle when it is
plugged in at time $\tau = A$. 

We can model electric vehicle charging as a deferrable load,
where we require a given charge level to be achieved at some time.
A more realistic model is as a combination of a deferrable and
curtailable load, with cost function 
\[
\alpha\sum_{\tau=A}^D(c^\mathrm{des}(\tau)-q(\tau))_+,
\]
where $\alpha > 0$ is a penalty parameter.
Here $c^\mathrm{des}(\tau)$ is the desired charge level at time $\tau$,
and $c^\mathrm{des}(\tau)-q(\tau))_+$ is the shortfall.

\paragraph{External tie with transaction cost.}
An external tie is a connection to an external source of power. We
represent this as a single terminal device with power schedule
$p_\mathrm{ex}$. 
In this case, $p_\mathrm{ex}(\tau)_- =
\max\{-p_\mathrm{ex}(\tau),0\}$ is the amount of energy pulled from the
source, and $p_\mathrm{ex}(\tau)_+ = \max\{p_\mathrm{ex}(\tau),0\}$ is
the amount of energy delivered to the source, at time $\tau$. We have
the constraint $|p_\mathrm{ex}(\tau)| \leq E^\mathrm{max}(\tau)$, where
$E^\mathrm{max}\in\reals^T$ is the transaction limit. 

We suppose that the prices for buying and selling energy are given by
$c \pm \gamma$ respectively, where $c(\tau)$ is the midpoint price, 
and $\gamma(\tau)>0$
is the difference between the price for buying and selling (\ie, the 
transaction cost).
The cost function is then
\[
-(c-\gamma)^T(p_\mathrm{ex})_++(c+\gamma)^T(p_\mathrm{ex})_- =
-c^Tp_\mathrm{ex}+\gamma^T|p_\mathrm{ex}|,
\]
where $|p_\mathrm{ex}|$, $(p_\mathrm{ex})_+$, and $(p_\mathrm{ex})_-$
are all interpreted elementwise.

\section{Convexity}\label{s-convexity}
\subsection{Devices}
We call a device convex if its objective function is a convex function.
A network is convex if all of its devices are convex. For convex
networks, the OPSP is a convex optimization problem, which means that in
principle we can efficiently find a global solution \cite{BoV:04}. When
the network is not convex, even finding a feasible solution for the OPSP
can become difficult, and finding and certifying a globally optimal
solution to the OPSP is generally intractable. However, special
structure in many practical power distribution problems allows us to
guarantee optimality in certain cases.

In the examples from \S \ref{s-devices}, the battery, fixed load,
thermal load, curtailable load, electric vehicle, and external tie are all
convex devices using the constraints and objective functions given.  A
deferrable load is convex if we drop the constraint that it can only be
turned on or off. We discuss the convexity properties of the generator
and transmission line in the following section.

\subsection{Relaxations}

One technique to deal with non-convex networks is to use convex
relaxations. We use the notation $g^\mathrm{env}$ to denote the convex
envelope \cite{Roc:70} of the function $g$. There are many equivalent
definitions for the convex envelope, for example, $g^\mathrm{env}$ =
$(g^*)^*$, where $g^*$ denotes the convex conjugate of the function $g$.
We can equivalently define $g^\mathrm{env}$ to be the largest convex
lower bound of $g$. If $g$ is a convex, closed, proper (CCP) function,
then $g=g^\mathrm{env}$.

We define the \emph{relaxed optimal power scheduling problem} (rOPSP) as
\begin{equation}\label{e-network-relaxed}
\begin{array}{ll}
\mbox{minimize} & f^\mathrm{env}(p)\\
\mbox{subject to} & \bar p = 0,
\end{array}
\end{equation}
with variable $p:\terminals \to \reals^T$. This is a convex optimization
problem, whose optimal value can in principle be computed efficiently,
and whose optimal objective value is a lower bound for the optimal
objective value of the OPSP. In some cases, we can guarantee \emph{a
priori} that a solution to the rOPSP will also be a solution to the
OPSP based on a property of the network objective such as monotonicity.
Even when the relaxed solution does not satisfy all of the constraints
in the unrelaxed problem, it can be used as a starting point to help
construct good, local solutions to the unrelaxed problem. The
suboptimality of these local solutions can then be bounded by the gap
between their network objective and the lower bound provided by the
solution to the rOPSP. If this gap is small for a given local solution,
we can guarantee that it is nearly optimal.

% In many practical cases, the solution to
% the rOPSP will satisfy some or all of the constraints in the original OPSP.

\paragraph{Generator.}

\begin{figure}
\begin{center}
\subfigure{
    \psfrag{000}[ct][ct]{$0$}%
    \psfrag{001}[ct][ct]{$1$}%
    \psfrag{002}[ct][ct]{$2$}%
    \psfrag{003}[ct][ct]{$3$}%
    \psfrag{004}[ct][ct]{$4$}%
    %\psfrag{005}[rc][rc]{$5$}%
    \psfrag{005}[rc][rc]{$0.4$}%
    \psfrag{006}[rc][rc]{$0.6$}%
    \psfrag{007}[rc][rc]{$0.8$}%
    \psfrag{008}[rc][rc]{$1.0$}%
    \psfrag{009}[rc][rc]{$1.2$}%
    \psfrag{010}[rc][rc]{$1.4$}%
    \psfrag{011}[rc][rc]{$1.6$}%
    \psfrag{pmin}[b][l]{$P^\mathrm{min}$}
    \psfrag{pmax}[][]{$P^\mathrm{max}$}
    \psfrag{pc}[][]{$P^\mathrm{c}$}
    \psfrag{pgen}[t][B]{$-p_\mathrm{gen}$}
    \psfrag{fgen}[b][t]{$\phi_\mathrm{gen}(p_\mathrm{gen})$}
    \includegraphics[width = 0.4\textwidth]{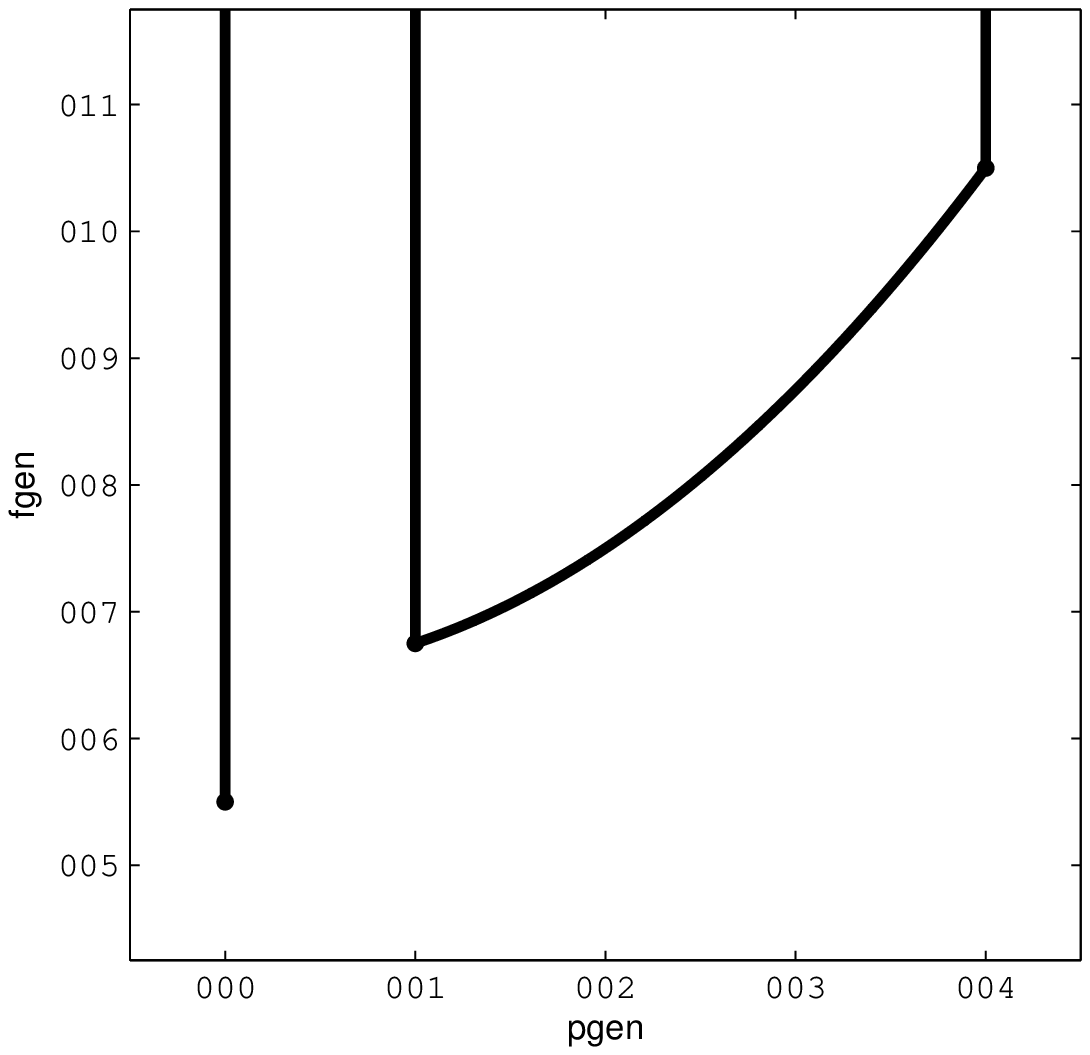}
            }
    \hspace{0.75cm}
\subfigure{
    \psfrag{000}[ct][ct]{$0$}%
    \psfrag{001}[ct][ct]{$1$}%
    \psfrag{002}[ct][ct]{$2$}%
    \psfrag{003}[ct][ct]{$3$}%
    \psfrag{004}[ct][ct]{$4$}%
    %\psfrag{005}[rc][rc]{$5$}%
    \psfrag{005}[rc][rc]{$0.4$}%
    \psfrag{006}[rc][rc]{$0.6$}%
    \psfrag{007}[rc][rc]{$0.8$}%
    \psfrag{008}[rc][rc]{$1.0$}%
    \psfrag{009}[rc][rc]{$1.2$}%
    \psfrag{010}[rc][rc]{$1.4$}%
    \psfrag{011}[rc][rc]{$1.6$}%
    \psfrag{pmin}[b][l]{$P^\mathrm{min}$}
    \psfrag{pmax}[][]{$P^\mathrm{max}$}
    \psfrag{pc}[][]{$P^\mathrm{c}$}
    \psfrag{pgen}[t][B]{$-p_\mathrm{gen}$}
    \psfrag{fgen}[b][t]{$\phi_\mathrm{gen}^\mathrm{env}(p_\mathrm{gen})$}
    \includegraphics[width = 0.4\textwidth]{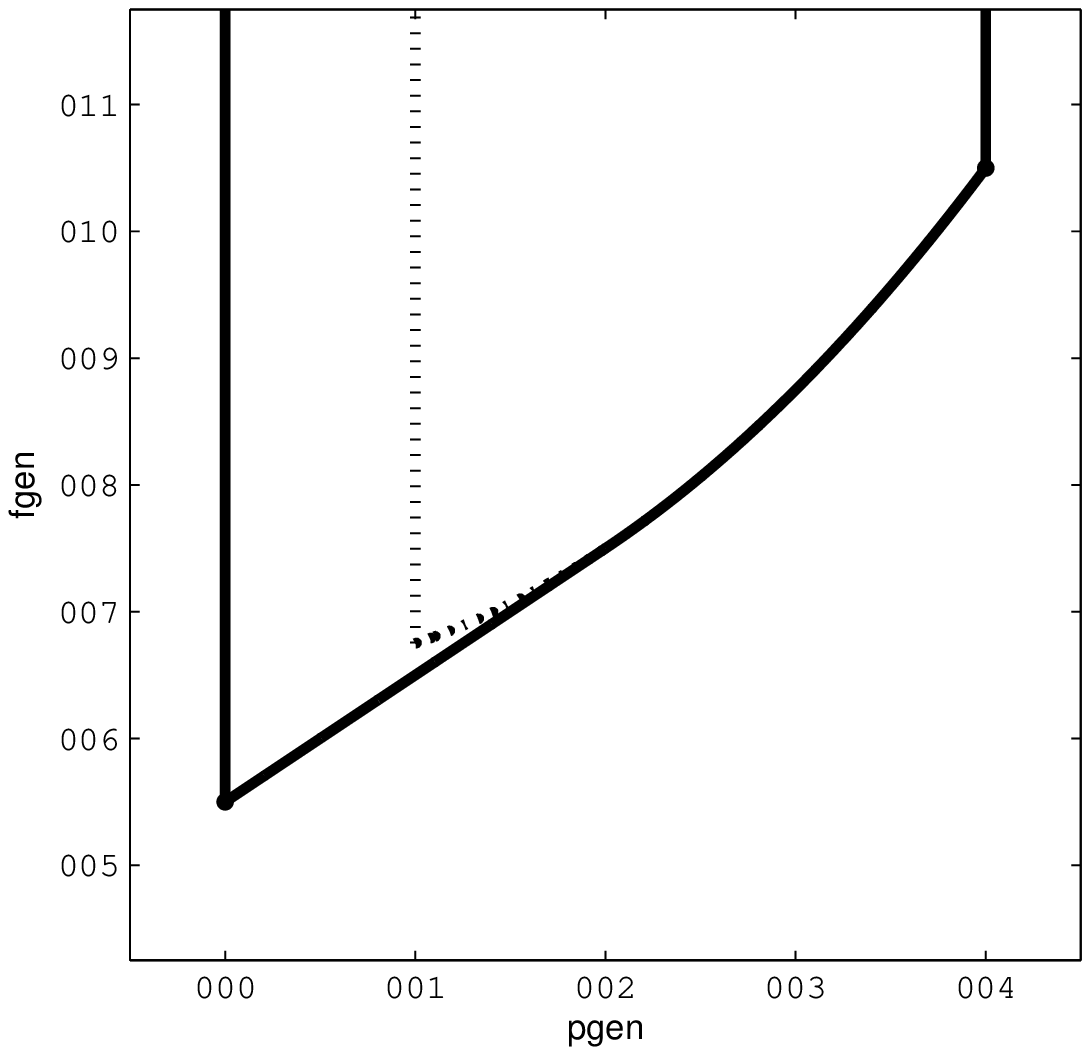}
            }
\end{center}
\caption{\emph{Left}: Cost function for a generator that can be turned
off. \emph{Right}: Its convex relaxation.}
\label{f-generator}
\end{figure}

When a generator is modeled as in \S \ref{s-devices} and is always
powered on, it is a convex device. However, when given the ability to be
powered on and off, the generator is no longer convex. In this case, we
can relax the generator objective function so that its cost for power
production in each time period, given in figure \ref{f-generator}, is a
convex function. This allows the generator to produce power in the
interval $[0,P^\mathrm{min}]$.

% The rOPSP can be used to solve the unit commitment problem (cite XXX) for
% generators.  In the solution to the rOPSP, if all generators have their power
% schedules either identically equal to zero or always in the range
% $[P^\mathrm{c},P^\mathrm{max}]$, and all other constraints in the original OPSP
% are satisfied, then the rOPSP gives a solution to the OPSP.  When this is not
% the case, the rOPSP solution can be used to partition the set of all generators
% into two sets: ON and OFF.  After this partitioning, a modified version of the
% OPSP is solved, which adds the constraints that all generators in the OFF set
% have $p_\mathrm{gen}=0$, while all generators in the ON set are constrained to
% operate in the range $P^\mathrm{min} \leq p_\mathrm{gen} \leq P^\mathrm{max}$.
% With these additional restrictions, all generators are convex devices, and
% (assuming no other non-convex devices), the resulting network is convex.
% Although solving this modified OPSP no longer is guaranteed to give the optimal
% value for the OPSP, it is in principle tractable to compute and can produce 
% good local solutions.

\paragraph{Transmission line.}
\begin{figure}
\begin{center}
\subfigure{
    \psfrag{x}[][]{}
    \psfrag{y}[][]{}
    \includegraphics[width = 0.4\textwidth]{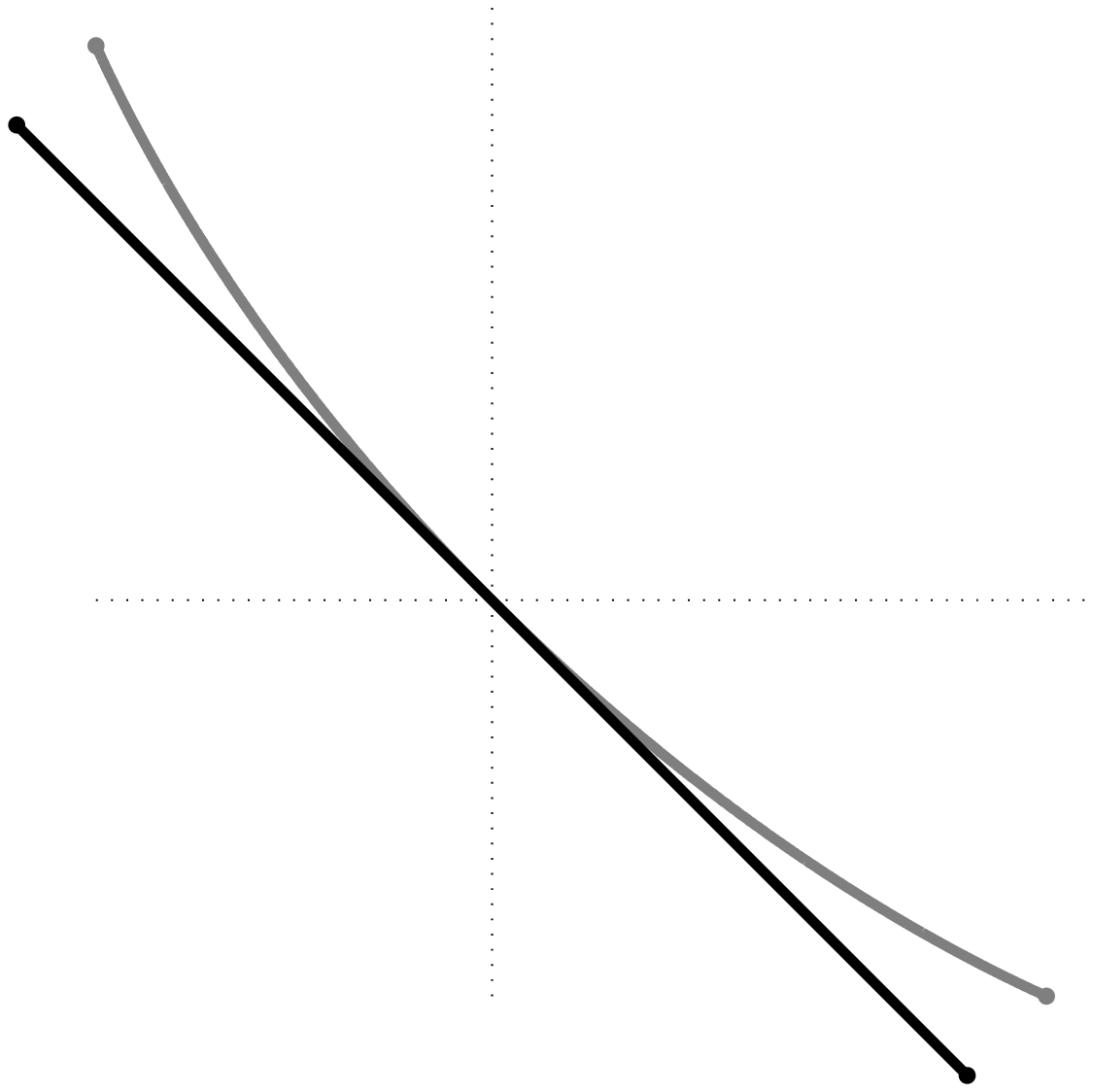}
    %XXX huge hack, fix later
    \rput(-0.8, 3.3){$p_1$}
    \rput(-3.3, 6.0){$p_2$}
}
\hspace{0.75cm}
\subfigure{
    \psfrag{x}[][]{}
    \psfrag{y}[][]{}
    \includegraphics[width = 0.4\textwidth]{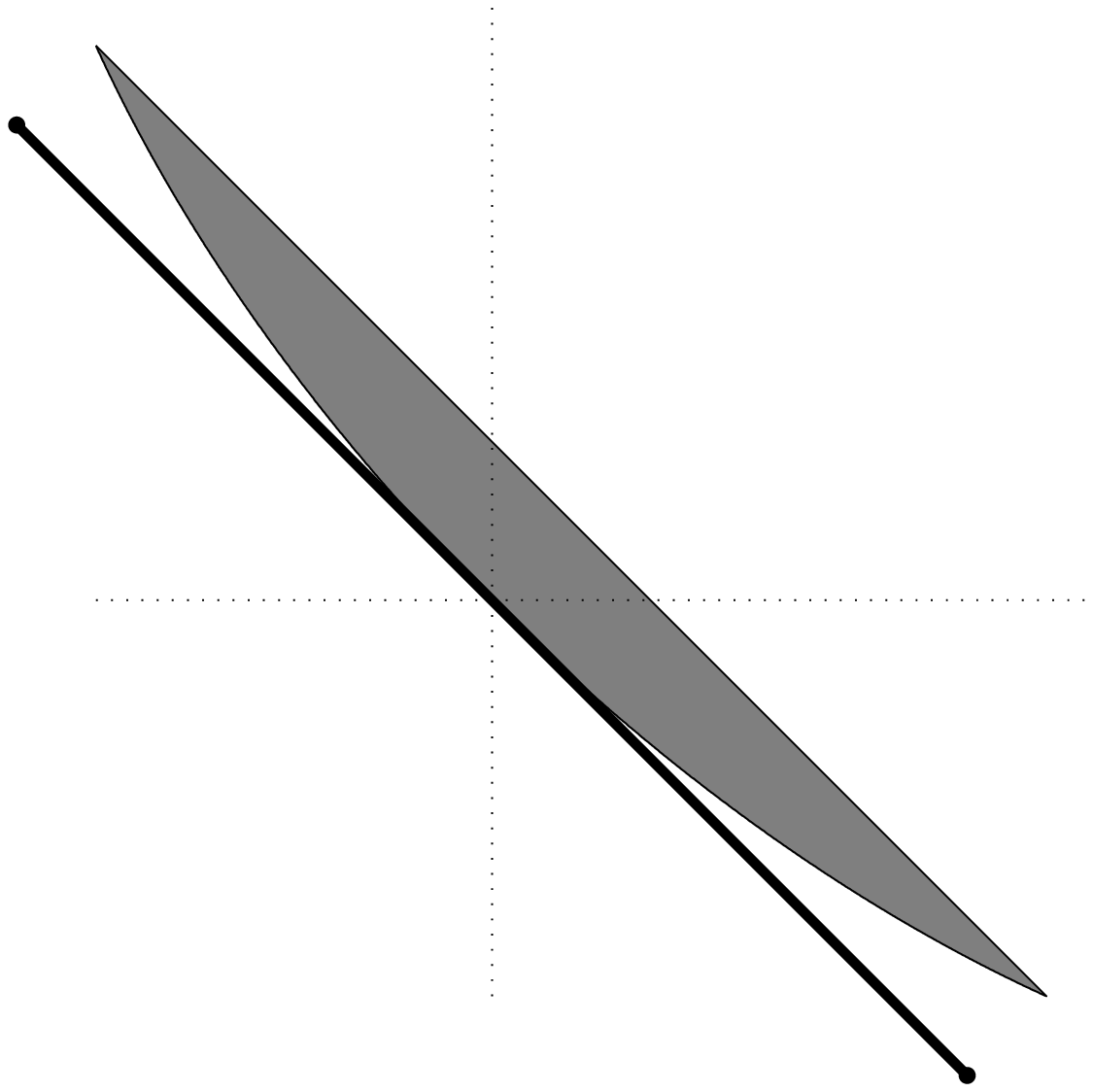}
    %XXX huge hack, fix later
    \rput(-0.8, 3.3){$p_1$}
    \rput(-3.3, 6.0){$p_2$}
}
\end{center}
\caption{\emph{Left}: Feasible sets of a transmission lines with no loss
(black) and AC loss (grey). \emph{Right}: Their convex relaxations.}
\label{f-line-loss}
\end{figure}

In a lossless transmission line, we have $\ell(p_1,p_2) = 0$, and thus
the set of feasible power schedules is the line segment
\[
L = \{ (p_1,p_2) \mid p_1=-p_2, \quad p_2 \in [-C^\mathrm{max}/2, C^\mathrm{max}/2]
\},
\]
as shown in figure \ref{f-line-loss} in black. When the transmission
line has losses, in most cases the loss function $\ell$ is a convex
function of the input and output powers, which leads to a feasible
region like the grey arc in figure \ref{f-line-loss}. For example, using
a lumped $\Pi$ model and under the common assumption that the voltage
magnitude is fixed \cite{BV:99}, a transmission line with series
admittance $g - ib$ gives the quadratic loss
\begin{equation} \label{e-line-loss}
\ell(p_1, p_2) = -(g/4)((p_1 + p_2)^2/g^2 + (p_1 -p_2)^2/b^2).
\end{equation}

The feasible set of a relaxed transmission line is given by the convex
hull of the original transmission line's constraints. The right side of
figure \ref{f-line-loss} shows examples of this for both lossless and
lossy transmission lines. Physically, this relaxation gives lossy
transmission lines the ability to discard some additional power beyond
what is simply lost to heat. Since electricity is generally a valuable
commodity in power networks, the transmission lines will generally not
throw away any additional power in the optimal solution to the rOPSP,
leading to the power line constraints in the rOPSP being tight and thus
also satisfying the unrelaxed power line constraints in the original
OPSP. As was shown in \cite{LTZ:12}, when the network is a tree, this
relaxation is always tight.  In addition, when all locational marginal prices
are positive and no other non-convexities exist in the network, the tightness
of the line constraints in the rOPSP can be guaranteed in the case of networks
that have separate phase shifters on each loop in the networks
whose shift parameter can be freely chosen \cite{SL:12a}.

\section{Decentralized method} \label{s-method}

We begin this section by deriving the prox-average message passing
equations assuming that all the device objective functions are convex
closed proper (CCP) functions. We then compare the computational and
communication requirements of prox-average message passing with a
centralized solver for the OPSP. The additional requirements that the
functions are closed and proper are technical conditions that are in
practice satisfied by any convex function used to model devices. We note
that we do not require either finiteness or strict convexity of any
device objective function, and that all results apply to networks with
arbitrary topology.

Whenever we have a set of variables that maps terminals to time periods,
$x:\terminals \to \reals^T$ (which we can also associate with a
$|\terminals| \times T$ matrix), we will use the same index and
over-line notation for the variables $x$ as we do for power schedules
$p$. For example, $x_t \in \reals^T$ consists of the time period vector
of values of $x$ associated with terminal $t$ and $\bar{x}_t =
(1/|n|)\sum_{t'\in n} x_{t'}$, where $t \in n$, with similar
notation for indexing $x$ by devices and nets.

\subsection{Prox-average message passing}

We derive the prox-average message passing equations by reformulating
the OPSP using the alternating direction method of multipliers (ADMM)
and then simplifying the resulting equations. We refer the reader to
\cite{BP:11} for a thorough overview of ADMM.

We first rewrite the OPSP as
\begin{equation}\label{e-opsp-split}
\begin{array}{ll}
\mbox{minimize} & \sum_{d\in\devices} f_d(p_d) +
		  \sum_{n\in\nets}    g_n(z_n) \\
\mbox{subject to} & p = z,
\end{array}
\end{equation}
with variables $p,z: \terminals \to \reals^T$, where $g_n(z_n)$ is the
indicator function on the set $\{z_n: \bar{z}_n = 0\}$. We use the
notation from \cite{BP:11} and, ignoring a constant, form the augmented
Lagrangian
\begin{eqnarray}\label{e-augmented-lagrangian}
L_\rho(p,z,u) = \sum_{d\in\devices} f_d(p_d) + \sum_{n\in\nets} g_n(z_n)
+ (\rho/2)\|p - z + u\|_2^2,
\end{eqnarray}
with the scaled dual variables $u = y/\rho : \terminals \to \reals^T$,
which we associate with a $|\terminals| \times T$ matrix. Because
devices and nets are each partitions of the terminals, the last term of
(\ref{e-augmented-lagrangian}) can be split across either devices
or nets, \ie,
\[
(\rho/2)\|p - z + u\|_2^2 =
\sum_{d\in\devices} (\rho/2)\|p_d - z_d + u_d\|_2^2 =
\sum_{n\in\nets} (\rho/2)\|p_n - z_n + u_n\|_2^2.
\]
The resulting ADMM algorithm is then given by the iterations
\begin{eqnarray*}
p_d^{k+1} &:=&  \argmin_{p_d} \left( f_d(p_d) + (\rho/2)
\|p_d-z_d^k + u_d^k\|^2_2\right), \quad d \in \devices, \\
z_n^{k+1} &:=& \argmin_{z_n} \left( g(z_n) +  (\rho/2)
\|z_n-u_n^k-p_n^{k+1}\|^2_2\right), \quad n \in \nets, \\
u_n^{k+1} &:=&  u_n^k + (p_n^{k+1}-z_n^{k+1}), \quad n \in \nets,
\end{eqnarray*}
where the first step is carried out in parallel by all devices, and then
the second and third steps are carried out in parallel by all nets.

Since $g_n(z_n)$ is simply an indicator function for each net $n$, the
second step of the algorithm can be computed analytically and is given
by
\begin{eqnarray*}
z^{k+1}_n &:=& u^k_n + p^{k+1}_n - \bar{u}^k_n - \bar{p}^{k+1}_n.
\end{eqnarray*}
By plugging this quantity into the $u$ update step, the algorithm can be
simplified further to yield the \textbf{prox-average message passing algorithm}:
\begin{enumerate}
\item \emph{Proximal power schedule update.}
\begin{equation}
p_d^{k+1} :=  \mathbf{prox}_{f_d, \rho}(p_d^k - \bar{p}_d^k - u_d^k),
\quad d\in \devices.
\end{equation}
\item \emph{Scaled price update.}
\begin{equation}
u_n^{k+1} :=  u_n^k + \bar p_n^{k+1}, \quad n \in \nets.
\end{equation}
\end{enumerate}
The proximal function for a function $f$ is given by
\begin{equation} \label{e-prox}
\mathbf{prox}_{f,\rho}(x) = \argmin_{y}
\left(f(y)+(\rho/2)\|x-y\|_2^2 \right),
\end{equation}
which is guaranteed to exist when $f$ is CCP \cite{Roc:70}.

We can now see where the name prox-average message passing comes from.
In each iteration, every device computes the proximal function of its
objective function, with an argument that depends on messages passed to
it through its terminals by its neighboring nets in the previous
iteration ($\bar{p}_d^k$ and $u_d^k$). Then, every devices passes to its
terminals the newly computed power schedules, $p_d^{k+1}$, which are
then passed to the terminals' associated nets. Every net computes its
new average power imbalance, $\bar{p}_n^{k+1}$, updates its dual
variables, $u_n^{k+1}$, and broadcasts these values to its associated
terminals' devices. Since $\bar{p}_n^k$ is simply $|n|$ copies of the same
vector for all $k$, we can see that all terminals connected to the same
net have the same value for their dual variables throughout the
algorithm, \ie, for all values of $k$, $u_t^k = u_{t'}^k$ whenever $t,
t' \in n$ for any $n \in \nets$.

As an example, consider the network represented by figures
\ref{f-network} and \ref{f-bipartite-graph}. The prox-average algorithm
performs the power schedule update on the devices (the boxes on the left
in figure \ref{f-bipartite-graph}). The devices share the respective
power profiles via the terminals, and the nets (the solid boxes on the
right) compute the scaled price update. For any network, the
prox-average algorithm can be thought of as alternating between devices
(on the left) and nets (on the right).

\paragraph{Convergence.}
We make a few comments about the convergence of prox-average message
passing. Since prox-average message passing is a version of ADMM, all
convergence results that apply to ADMM also apply to prox-average
message passing. In particular, when all devices have closed, convex,
proper (CCP) objective functions and a feasible solution to the OPSP
exists, the following hold.
\begin{enumerate}
\item \emph{Residual convergence.}
$\bar{p}^k \to 0$ as $k\to\infty$,
\item \emph{Objective convergence.}
$\sum_{d\in\devices} f_d(p_d^k)+\sum_{n\in\nets}g_n(p_n^k)\to f^\star$ 
as $k\to\infty$,
\item \emph{Dual variable convergence.}
$\rho u^k = y^k \to y^\star$ as $k\to\infty$,
\end{enumerate}
where $f^\star$ is the optimal value for the OPSP, and $y^\star$ are the
optimal dual variables. The proof of these conditions can be found in
\cite{BP:11}. As a result of the third condition, the optimal locational
marginal prices $\mathcal{L}^\star$ can be found for each net
$n\in\nets$ by setting $\mathcal{L}_n^\star = y_n^\star/|n|$.

\paragraph{Stopping criterion.}
Following \cite{BP:11}, we can define primal and dual residuals,
which for prox-average message passing simplify to
\[
r^k = \bar{p}^k, \quad 
s^k = \rho((p^k-\bar{p}^k)-(p^{k-1}-\bar{p}^{k-1})).
\]
We give a simple interpretation of each residual. The primal residual is
simply the net power imbalance across all nets in the network, which is
the original measure of primal feasibility in the OPSP. The dual
residual is equal to the difference between the current and previous
iterations of the difference between power schedules and their average
net power. The locational marginal price at each net is determined by
the deviation of all associated terminals' power schedule from the
average power on that net. As the change in these deviations approaches
zero, the corresponding locational marginal prices converge to their
optimal values.

We can define a simple criterion for terminating prox-average message
passing when
\[
\|r^k\|_2 \leq \epsilon^\mathrm{pri},\quad 
\|s^k\|_2 \leq \epsilon^\mathrm{dual},
\]
where $\epsilon^\mathrm{pri}$ and $\epsilon^\mathrm{dual}$ are,
respectively, primal and dual tolerances. We can normalize both of these
quantities to network size by the relation
\[
\epsilon^\mathrm{pri} = 
\epsilon^\mathrm{dual} = 
\epsilon^\mathrm{abs}\sqrt{|\terminals|T},
\]
for some absolute tolerance $\epsilon^\mathrm{abs} > 0$.

% XXX - relative stopping criteria by broadcasting local terminal primal and dual
% residuals, along with the number of hops the residual has gone through
% Thus, if not message is above the a certain residual threshold and is newer
% than say twice the diameter (eg 20) number of hops, we can terminate.
% In practice, the norms of the primal and dual residuals can be spread
% throughout the network using gossip algorithms, or individual devices can
% simply use local residuals in determining when to terminate their message
% passing.

\paragraph{Choosing a value of $\rho$.}
Numerous examples show that the value of $\rho$ can have a dramatic
effect on the rate of convergence of ADMM and prox-average message
passing. Many good methods for picking $\rho$ in both offline and online
fashions are discussed in \cite{BP:11}. We note that unlike other
versions of ADMM, the scaling parameter $\rho$ enters very simply into
the prox-average equations and can thus be modified online without
incurring any additional computational penalties, such as having to
re-factorize a matrix.

We can modify the prox-average message passing algorithm with the
addition of a third step
\begin{enumerate}
\setcounter{enumi}{2}
\item \emph{Parameter update and price rescaling.}
\begin{eqnarray*}
\rho^{k+1} &:=& h(\rho^k, r^k, s^k), \\
u^{k+1} &:=& \frac{\rho^k}{\rho^{k+1}} u^{k+1},
\end{eqnarray*}
\end{enumerate}
for some function $h$.
We desire to pick an $h$ such that the primal and dual residuals are of
similar size throughout the algorithm, \ie, $\rho^k \|r^k\|_2 \approx
\|s^k\|_2$ for all $k$. To accomplish this task, we use a simple
proportional-derivative controller to update $\rho$, choosing $h$ to be
\[
h(\rho^k) = \rho^k \exp(\lambda v^k + \mu(v^k - v^{k-1})),  
\]
where $v^k = \rho^k \|r^k\|_2/\|s^k\|_2 - 1$ and $\lambda$ and $\mu$ are
nonnegative parameters chosen to control the rate of convergence.
Typical values of $\lambda$ and $\mu$ are between $10^{-3}$ and
$10^{-1}$.

When $\rho$ is updated in such a manner, convergence is sped up in many
examples, sometimes dramatically. Although it can be difficult to prove
convergence of the resulting algorithm, a standard trick is to assume
that $\rho$ is changed only for a large but bounded number of
iterations, at which point it is held constant for the remainder of the
algorithm, thus guaranteeing convergence.

% (XXX: I also clip $\rho^k$ to be between $\epsilon^\mathrm{abs}$ and
% $1/\epsilon^\mathrm{abs}$---not sure if that should be mentioned?)

\paragraph{Non-convex case.}

When one or more of the device objective functions is non-convex, we can
no longer guarantee that prox-average message passing converges to the
optimal value of the OPSP or even that it converges at all (\ie, reaches
a fixed point). Prox functions for non-convex devices must be carefully
defined as the set of minimizers in (\ref{e-prox}) is no longer
necessarily a singleton. Even when they can be defined, prox functions
of non-convex functions are intractable to compute in many cases.

One solution to these issues is to use prox-average message passing to
solve the rOPSP. It is easy to show that $f^\mathrm{env}(p) = \sum_{d
\in \devices} f^\mathrm{env}_d(p_d)$. As a result, we can run
prox-average message passing using the proximal functions of the relaxed
device objective functions. Since $f_d^\mathrm{env}$ is a CCP function
for all $d\in\devices$, prox-average message passing in this case is
guaranteed to converge to the optimal value of the rOPSP and yield the
optimal relaxed locational marginal prices.

\subsection{Discussion}

In order to compute the prox-average messages, devices and nets only
require knowledge of who their network neighbors are, the ability to
send small vectors of numbers to those neighbors in each iteration, and the
ability to store small amounts of state information and efficiently
compute prox functions (devices) or averages (nets). As all
communication is local and peer-to-peer, prox-average message passing
supports the ad hoc formation of power networks, such as micro grids,
and is robust to device failure and unexpected network topology changes.

Due to recent advances in convex optimization \cite{WB:10, CVXGEN, MWB:11},
many of the prox function calculations that devices must perform can be
very efficiently executed at millisecond or microsecond time-scales on
inexpensive, embedded processors. Since all devices and all nets can
each perform their computations in parallel, the time to execute a
single, network wide prox-average message passing iteration (ignoring
communication overhead) is equal to the sum of the maximum computation
time over all devices and the maximum computation time of all nets
in the network. As a result, the computation time per iteration is small
and essentially independent of the size of the network.

In contrast, solving the OPSP in a centralized fashion requires complete
knowledge of the network topology, sufficient communication bandwidth to
centrally aggregate all devices objective function data, and sufficient
centralized computational resources to solve the resulting OPSP. In
large, real-world networks, such as the smart grid, all three of these
requirements are generally unattainable. Having accurate and timely
information on the global connectivity of all devices is infeasible for
all but the smallest of dynamic networks. Centrally aggregating all
device objective functions would require not only infeasible bandwidth
and data storage requirements at the aggregation site, but also the
willingness of all devices to expose what could be proprietary function
parameters in their objective functions. Finally, a centralized solution to the
OPSP requires solving an optimization problem with
$\Omega(|\terminals|T)$ variables, which leads to an identical lower
bound on the time scaling for a centralized solver, even if problem
structure is exploited.  As a result, the centralized solver cannot scale to
solve the OPSP on very large networks.

\section{Numerical examples} \label{s-numerical}

We illustrate the speed and scaling of prox-average message passing with
a range of numerical examples. In the first two sections, we describe
how we generate network instances for our examples. We then describe our
implementation, showing how multithreading can exploit problem
parallelism and how our method would scale in a fully peer-to-peer
implementation. Lastly, we present our results, and demonstrate how the
number of prox-average iterations needed for convergence is essentially
independent of network size and also significantly decreases when the
algorithm is seeded with a reasonable warm-start.

% We first solve the network energy management problem in a centralized
% fashion, and then show that our decentralized algorithm quickly
% converges to the same solution. In addition, we show that the
% networked devices quickly and automatically adapt to dynamic changes to the
% network, including the addition of new devices, and both planned and 
% unexpected topology changes.

\subsection{Network topology}
We generate a network instance by first picking the number of nets $N$.
We generate the locations $x_i \in
\reals^2$, $i=1,\ldots, N$ by drawing them uniformly at random from $[0,
\sqrt{N}]^2$. 
(These locations will be used to determine network topology.)
Next, we introduce transmission lines into the network as
follows. We first connect a transmission line between all pairs of nets
$i$ and $j$ independently and with probability
\[
\gamma(i,j) = \alpha \min(1,d^2/\|x_i - x_j\|_2^2).
\]
In this way, when the distance between $i$ and $j$ is smaller than $d$,
they are connected with a fixed probability $\alpha > 0$, and when they
are located farther than distance $d$ apart, the probability decays as
$1/\|x_i-x_j\|_2^2$. After this process, we add a transmission line
between any isolated net and its nearest neighbor. We then introduce
transmission lines between distinct connected components by selecting
two connected components uniformly at random and then selecting two
nets, one inside each component, uniformly at random and connecting them
by a transmission line. We continue this process until the network is
connected.

% with probability $p^\mathrm{device}$, we
% connect either a generator, battery, fixed load, deferrable load, or
% curtailable load independently and uniformly at random. If a generator is
% chosen, we select its type (small, medium, or large) independently and
% uniformly at random as well.

For the network instances we present, we chose parameter values $d =
0.15$ and $\alpha = 0.8$ as the parameters for generating our network.
This results in networks with an average degree of $2.3$. Using these
parameters, we generated networks with $100$ to $100000$ nets, which
resulted in optimization problems with approximately $30$ thousand to
$30$ million variables.

\begin{figure} \centering
    \includegraphics[width=0.4\textwidth]{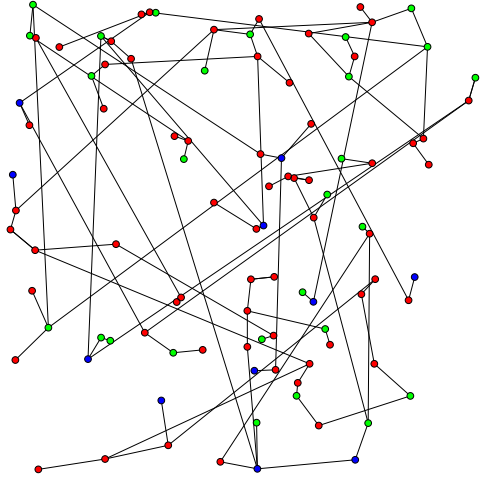}
    \caption{A sample random network. Devices are color-coded: 
	generators are in green, batteries are in blue, and loads are in 
	red. Edges represent transmission lines.}
    \label{f-network-example}
\end{figure}

\subsection{Devices}
After we generate the network topology described above, we randomly
attach a single (one-terminal) device to each net according to the
distribution in table \ref{t-device-prob}. The models for each device in
the network are identical to the ones given in section \ref{s-devices},
with model parameters chosen in a manner we describe below. Figure
\ref{f-network-example} shows an example network for $N = 100$ ($30$
thousand variables) generated in this fashion.

For simplicity, our examples only include networks with the devices
listed below. For all devices, the time horizon was chosen to be $T=96$,
indicating $15$ minute intervals for a $24$ hour power schedule, with
the time period $\tau=1$ corresponding to midnight.

\begin{table}
	\centering
	\begin{tabular}{ l | c }
	    \textbf{Device} & \textbf{Fraction} \\ \hline \hline 
		Generator & $0.2$ \\
		Battery & $0.1$ \\ 
		Fixed load & $0.5$ \\ 
		Deferrable load & $0.1$ \\
		Curtailable load & $0.1$ \\
	\end{tabular}
	\caption{Fraction of one-terminal devices present in the generated networks.}
	\label{t-device-prob}
\end{table}

\paragraph{Generator.} Generators have the quadratic cost functions
given in \S \ref{s-devices} and are divided into
three types: small, medium, and large.
% and constraints with the same form as given in \S
% \ref{s-devices}, but with different values for the parameters $P^\mathrm{min},
% P^\mathrm{max}, R^\mathrm{max}, \alpha$, and $\beta$.
In each case, we allow the generator to be turned on and off by setting
$P_\mathrm{min} = 0$. Small generators have the smallest maximum power
output, but the largest ramp rates, while large generators have the
largest maximum power output, but the slowest ramp rates. Medium
generators lie in between. Large generators are generally more efficient
than small and medium generators which is reflected in their cost
function by having smaller values of $\alpha$ and $\beta$. Whenever a
generator is placed into a network, its type is selected uniformly at
random, and its parameters are taken from the appropriate row in table
\ref{t-generator}.

\begin{table}
	\centering
	\begin{tabular}{ l | c | c | c | c | c }
		& $P^\mathrm{min}$ & $P^\mathrm{max}$ & $R^\mathrm{max}$ & $\alpha$ & $\beta$ \\
		\hline\hline
		Large & $0$ & $50$ & $3$ & $0.001$ & $0.1$ \\
		Medium & $0$ & $20$ & $5$ & $0.005$ & $0.2$\\
		Small & $0$ & $10$ & $10$ & $0.02$ & $1$
	\end{tabular}
	\caption{Generator parameters.}
	\label{t-generator}
\end{table}

\paragraph{Battery.} For a given instance of a battery, its parameters
are generated by setting $q^\mathrm{init} = 0$ and selecting
$Q^\mathrm{max}$ uniformly at random from the interval $[20,50]$. The
charging and discharging rates are selected to be equal (\ie,
$C^\mathrm{max} = D^\mathrm{max}$) and drawn uniformly at random from
the interval $[5,10]$.

\paragraph{Fixed load.} The load profile for a fixed load instance is a
sinusoid,
\[
l(\tau) = c + a\sin(2\pi(\tau - \phi_0)/T), \quad \tau=1, \ldots, T,
\] 
with the amplitude $a$ chosen uniformly at random from the interval
$[1,5]$, and the DC term $c$ chosen so that $c = a + u$, where $u$ is
chosen uniformly at random from the interval $[0,0.5]$, which ensures
that the load profile remains elementwise positive. The phase shift
$\phi_0$ is chosen uniformly at random from the interval $[60,72]$,
ensuring that the load profile peaks between the hours of $3$pm and
$6$pm.

\paragraph{Deferrable load.} For an instance of a deferrable load, we
choose $E$ uniformly at random from the interval $[500,1000]$. The start
time index $A$ is chosen uniformly at random from the discrete set
$\{1,\ldots, T-7\}$. The end time index $D$ is then chosen uniformly at
random over the set $\{A+7,\ldots,T\}$, so that the minimum time window
to satisfy the load is $8$ time periods ($2$ hours). We set the maximum
power so that it is possible to satisfy the total energy constraint by only
operating in half of the available time periods, \ie,
$L^\mathrm{max} = 2E/(D-A)$.

\paragraph{Curtailable loads.} For an instance of a curtailable load,
the desired load $l$ is constant over all time periods with a magnitude
chosen uniformly at random from the interval $[5, 15]$. The penalty
parameter $\alpha$ is chosen uniformly at random from the interval
$[1,2]$.

\paragraph{Transmission line.} For an instance of a transmission line,
we choose its parameters by first solving the OPSP with lossless,
uncapacitated lines, where we add a small quadratic cost function
$\epsilon(p_1^2 + p_2^2)$, with $\epsilon = 10^{-3}$, to each
transmission line in order to help spread power flow throughout the
network. Using the flow values given by the solution to that problem, we
set $C^\mathrm{max} = \max(10, 4 F^\mathrm{max})$ for each line, where
$F^\mathrm{max}$ is equal to the maximum flow (from the solution to the
lossless problem) along that line over all time periods. We use the loss
function for transmission lines with a series admittance $g - ib$ given
by (\ref{e-line-loss}). We choose $g$ and $b$ such that $b = \gamma g$,
where $\gamma$ is chosen uniformly at random from the interval $[4.5,
5.5]$; on average, the susceptance is $5$ times larger than the
conductance. After we pick $C^\mathrm{max}$ and $\gamma$, the values of
$g$ and $b$ are chosen such that the loss when transmitting power at
maximum capacity is uniformly at random between $5$ to $15$ percent of
$C^\mathrm{max}$.

\subsection{Serial multithreaded implementation}
Our OPSP solver is implemented in C++, with the core prox-average
message passing equations occupying fewer than $25$ lines of C++
(excluding problem setup and classes). The code is compiled with
\verb|gcc| $4.4.5$ on an $8$-core, $3.4$GHz Intel Xeon processor with
$16$GB of RAM running the Debian OS. We used the compiler option
\verb|-O3| to leverage full code optimization.

To approximate a fully distributed implementation, we use \verb|gcc|'s
implementation of OpenMP (version $3.0$) and multithreading to
parallelize the computation of the prox functions for the devices. We
use $8$ threads (one per core) to solve each example. Assuming perfect
load balancing, this means that $8$ prox functions are being evaluated
in parallel. Effectively, we evaluate the prox functions by stepping
serially through the devices in blocks of size $8$. We do \emph{not},
however, parallelize the computation of the dual update over the nets
since the overhead of spawning threads dominates the vector operation
itself.

The prox functions for fixed loads and curtailable loads are separable
over $\tau$ and can be computed analytically. For more complex devices,
such as a generator, battery, or deferrable load, we compute the
prox function using CVXGEN \cite{CVXGEN}. The prox function for a
transmission line is computed by projecting onto the convex hull of the
line constraints.

For a given network, we solve the associated OPSP with an absolute
tolerance $\epsilon^\mathrm{abs} = 10^{-3}$. This translates to three
digits of accuracy in the solution. The CVXGEN solvers used to evaluate
the prox operators for some devices have an absolute tolerance of
$10^{-8}$. For our $\rho$-update function, $h$, we use the parameter
values $\lambda = 0.005$ and $\mu = 0.01$ and clip our values of $\rho$
to be between $\epsilon^\mathrm{abs}$ and $1/\epsilon^\mathrm{abs}$ to prevent
roundoff error.

\subsection{Peer-to-peer implementation} We have not yet created a fully
peer-to-peer, bulk synchronous parallel \cite{Val:90, MAB:10}
implementation of prox-average message passing, but have carefully
tracked prox-average solve times in our serial implementation in order
to facilitate a first order analysis of such a system. In a peer-to-peer
implementation, the proximal power schedule updates occur in parallel
across all devices followed by (scaled) price updates occurring in
parallel across all nets. As previously mentioned, the computation time
per prox-average iteration is thus the maximum time, over all devices,
to evaluate the proximal function of their objective, added to the
maximum time across all nets to average their terminal power schedules
and add that quantity to their existing price vector. Since evaluating
the prox function for some devices requires solving a convex
optimization problem, whereas the price update only requires a small
number of vector operations that can be performed as a handful of SIMD
instructions, the compute time for the price update is negligible in
comparison to the proximal power schedule update. The determining factor
in solve time, then, is in evaluating the prox functions for the power
schedule update. In our examples, the \emph{maximum} time taken to
evaluate any prox function is $1$ ms. To solve a problem with $N=100000$
nets ($30$ million variables) with approximately $500$ iterations of our
prox-average algorithm then takes only $500$ ms.

In practice, the actual solve time would clearly be dominated by network
communication latencies and actual runtime performance will be
determined by how quickly and reliably packets can be delivered. As a
result, in a true peer-to-peer implementation, a negligible amount of
time is actually spent on computation. However, it goes without saying
that many other issues must be addressed with a peer-to-peer protocol,
including handling network delays and security.

\subsection{Results}

\begin{figure}
	\begin{center}
	\subfigure{
		\centering
		\psfrag{iter}[t][b]{\makebox(0,20)[1]{\small iter $k$}}
		\psfrag{obj}[b][t]{\makebox(0,50)[1]{\small $|f^k - f^\star|/f^\star$}}
		\psfrag{1}[r][r]{\small $10^{-7}$}
		\psfrag{2}[r][r]{\small $10^{-5}$}
		\psfrag{3}[r][r]{\small $10^{-3}$}
		\psfrag{4}[r][r]{\small $10^{-2}$}
		\psfrag{5}[r][r]{\small $10^1$}
		\psfrag{0}[t][b]{\small $0$}
		\psfrag{250}[t][b]{\small $250$}
		\psfrag{500}[t][b]{\small $500$}
		\psfrag{750}[t][b]{\small $750$}
		\psfrag{1000}[t][b]{\small $1000$}
		\includegraphics[width=0.4\textwidth]{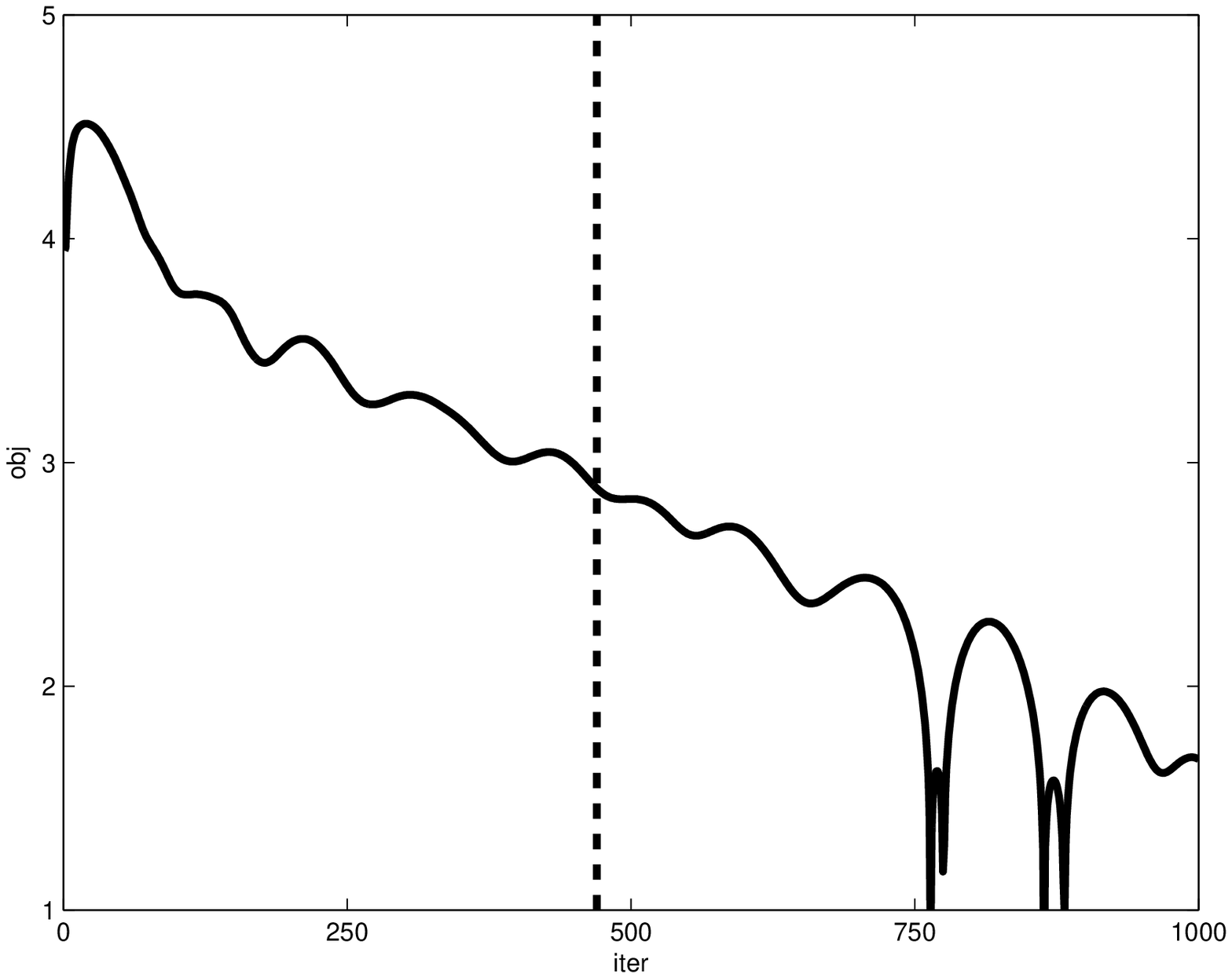}
	}
	\hspace{0.75cm}
	\subfigure{
		\centering
		\psfrag{iter}[t][b]{\makebox(0,20)[1]{\small iter $k$}}
		\psfrag{pbar}[b][t]{\makebox(0,50)[1]{\small $\| \overline p^k \|_2/\sqrt{|\mathcal T|T}$}}
		\psfrag{1}[r][r]{\small $10^{-7}$}
		\psfrag{2}[r][r]{\small $10^{-5}$}
		\psfrag{3}[r][r]{\small $10^{-3}$}
		\psfrag{4}[r][r]{\small $10^{-2}$}
		\psfrag{5}[r][r]{\small $10^1$}
		\psfrag{0}[t][b]{\small $0$}
		\psfrag{250}[t][b]{\small $250$}
		\psfrag{500}[t][b]{\small $500$}
		\psfrag{750}[t][b]{\small $750$}
		\psfrag{1000}[t][b]{\small $1000$}
		\includegraphics[width=0.4\textwidth]{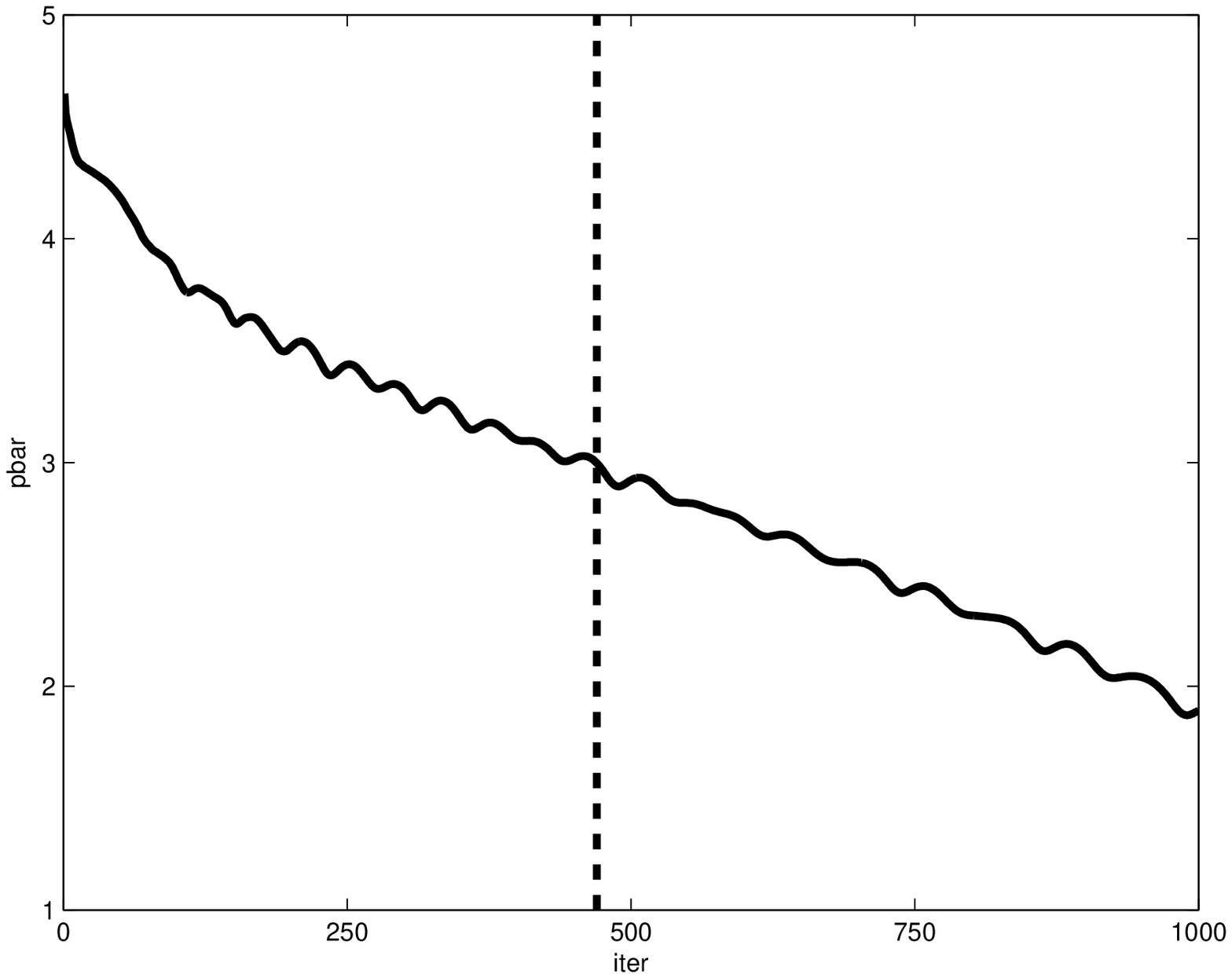}
	}
	\end{center}
\caption{The relative suboptimality (\emph{left}) and primal
infeasibility (\emph{right}) of prox-average message passing on a
network instance with $N=3000$ nets ($1$ million variables). The dashed
line shows when the stopping criterion is satisfied.}
\label{f-convergence}
\end{figure}

We first consider a single example: a network instance with $N=3000$
($1$ million variables). Figure \ref{f-convergence} shows that after
fewer than $500$ iterations of prox-average message passing, both the
relative suboptimality and the average net power imbalance are both less
than $10^{-3}$. The convergence rates for other network instances over
the range of sizes we simulated are similar.

\begin{figure}
	\centering
	\psfrag{N}[t][b]{$N$}
	\psfrag{time}[b][t]{time (seconds)}
	\psfrag{1}[r][r]{\small $1$}
	\psfrag{10}[r][r]{\small $10$}
	\psfrag{100}[r][r]{\small $100$}
	\psfrag{1000}[r][r]{\small $1000$}
	\psfrag{10000}[r][r]{\small $10000$}
	\psfrag{20}[t][]{\small $100$}
	\psfrag{200}[t][]{\small $1000$}
	\psfrag{2000}[t][]{\small $10000$}
	\psfrag{20000}[t][]{\small $100000$}
	\includegraphics[width=0.6\textwidth]{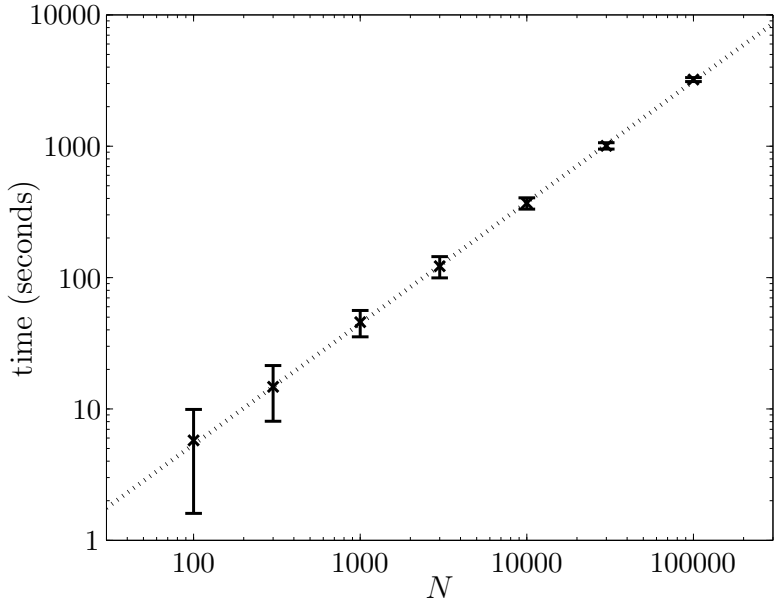}
    \caption{Average prox-average execution times for a family of 
	networks on $8$ cores. Error bars show $95\%$ confidence bounds. The 
	dotted line shows the least-squares fit to the data on a log-log 
	scale, resulting in an exponent of $0.923$.}
    \label{f-runtime}
\end{figure}

In figure \ref{f-runtime}, we present average timing results for solving
the OPSP for a family of examples, with networks of size $N=100$, $300$,
$1000$, $3000$, $10000$, $30000$, and $100000$. For each network size,
we generated and solved $20$ network instances to generate average solve
times and confidence intervals around those averages. For network
instances with $N=100000$ nets, the problem has approximately $30$
million variables, which we solve serially using prox-average message
passing in $52$ minutes on average.

For a peer-to-peer implementation, the runtime of prox-average message
passing should be essentially constant, and in particular independent of
the size of the network. For our multithreaded implementation, with
bounded computation, this would be reflected by a runtime that only
increases linearly with the number of nets in a network instance. By
fitting a line to figure \ref{f-runtime}, we find that our parallel
implementation scales as $O(N^{0.923})$. The slight discrepancy between
this and the ideal exponent of $1$ is accounted for by implementation
details such as operating system background processes consuming some
compute cycles and slightly imperfect load balancing across all $8$
cores in our system, especially for smaller values of $N$.

\begin{figure}
	\centering
	\psfrag{iters}[b][t]{\makebox(0,40)[1]{\small $K^\mathrm{warm}/K^\mathrm{cold}$}}
    \psfrag{sigma}[t][b]{\makebox(0,20)[1]{\small $\sigma$}}
	\psfrag{0}[r][]{\small $0$}
	\psfrag{0.2}[r][]{\small $0.2$}
	\psfrag{0.4}[r][]{\small $0.4$}
	\psfrag{0.6}[r][]{\small $0.6$}
	\psfrag{0.8}[r][]{\small $0.8$}
	\psfrag{1}[r][]{\small $1.0$}
	\psfrag{10}[t][]{\small $0.00$}
	\psfrag{20}[t][]{\small $0.05$}
	\psfrag{30}[t][]{\small $0.10$}
	\psfrag{40}[t][]{\small $0.15$}
	\psfrag{50}[t][]{\small $0.20$}
	%\psfrag{60}[t][]{\small $0.25$}

	\includegraphics[width=0.4\textwidth]{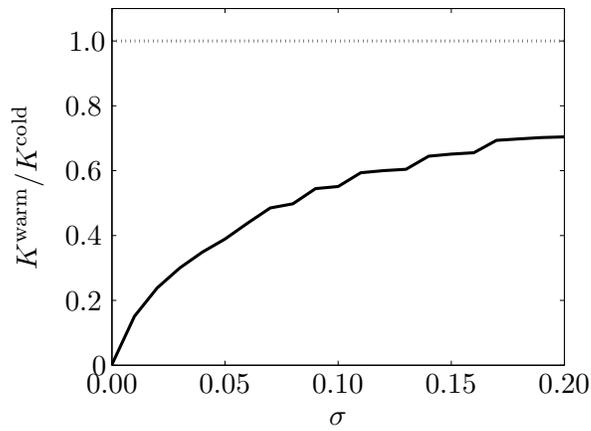}
	\vspace{0.5cm}
    \caption{Number of warm start iterations needed to converge for various perturbations of load profiles.}
    \label{f-warmstart}
\end{figure}

We note that figure \ref{f-runtime} shows \emph{cold start} runtimes for
solving the OPSP. If we have access to a good estimate of the power
schedules and locational marginal prices for each terminal, we can use
these estimates to \emph{warm start} our OPSP solver. To show the effect
of warm-starting, we solve a specific problem instance with $N=3000$
nets ($1$ million variables). We define $K^\mathrm{cold}$ to be the
number of iterations needed to solve an instance of this problem. We
then uniformly scale the load profiles of each device by separate and
independent lognormal random variables. The new load profiles, $\hat l$,
are obtained from the original load profiles $l$ according to
\[
\hat l = l \exp(\sigma X),
\]
where $X \sim \mathcal{N}(0,1)$, and $\sigma>0$ is given. Using the
solution of the original problem to warm start our solver, we solve the
perturbed problem and report the number of iterations $K^\mathrm{warm}$
needed to solve it. Figure \ref{f-warmstart} shows the ratio
$K^\mathrm{warm}/K^\mathrm{cold}$ as we vary $\sigma$, and indicates the
significant computational savings that warm-starting can achieve, even under
relatively large perturbations.

% We expect similar performance for other reasonable warm start estimates.

% In figure \ref{f-warmstart}, we uniformly
% scale the load profiles of each device by separate and independent lognormal
% random variables in a network with $N=3000$ nets and warm start the OPSP solver
% with the original solution.

\section{Extensions}
\label{s-extensions}

% \paragraph{AC power flow and OPF.}
% While we have only considered simple, scalar power flow for the examples
% in this paper, prox-average message passing can easily be extended to
% cover AC power as well. Under these circumstances, the feasible regions
% for transmission lines are ellipses and the OPSP can be interpreted as a
% dynamic version of the OPF over a time horizon. The results in
% \cite{SL:12a, LTZ:12} give broad conditions where solving the relaxed
% version of the OPF results in an exact solution for the OPF, which can
% be adapted to show that solving the rOPSP is equivalent to solving the
% OPSP under the similar conditions.
% 
% (XXX: but we solved this. we just didn't have reactive components.)

% As a result, semidefinite
% programming is no longer required to solve the OPF, so individual devices with
% a large number of states can ...

\paragraph{Receding horizon control.}
The speed with which prox-average message passing converges on very
large networks shows its applicability in coordinating real time
decisions across massive networks of devices. A direct extension of our
work to real-time network operation can be achieved using receding
horizon control (RHC) \cite{JMM:02, Bem:06}. In RHC, we solve the OPSP
at each time step to determine conditional consumption and generation
profiles for every device over the next $T$ time periods. We then
execute the first step of these profiles and resolve the OPSP using new
information and measurements that have become available. RHC has been
successfully applied to a wide range of areas, including power systems,
and allows us to take advantage of warm starting our algorithm, which we
have shown to significantly decrease the number of iterations needed for
convergence.

\paragraph{Hierarchical models.}
The power gird has a natural hierarchy, with generation and transmission
occurring at the highest level and residential consumption and
distribution occurring at the most granular. Prox-average message
passing can be easily extended into hierarchical interactions by
scheduling messages on different time scales and between systems at the
similar levels of the hierarchy \cite{CLCD:07}.

We can recursively apply prox-average message passing at each level of
the hierarchy. At the highest level, all regional systems exchange their
proximal updates once they have computed their own prox-function. It can
be shown that computing this function for a given region is equivalent
to computing a partial minimization over the sum of the objective
functions of devices located inside that region, subject to intra-region
power balance. This too can be computed using prox-average message
passing. This process can be continued down to the individual device
level, at which point the device must solve its own prox function
directly as the base case.

\paragraph{Local stopping criteria and $\rho$ updates.}
The stopping criterion and the algorithm we propose to update $\rho$ in
\S \ref{s-method} both currently require global device coordination
--- specifically the global values of the primal and dual residuals at
each iteration. In principle, these could be computed in a decentralized
fashion across the network by gossip algorithms \cite{Shah:08}, but
that would require many rounds of gossip in between each iteration of
prox-average message passing, significantly increasing runtime. We are
currently investigating methods by which both the stopping criterion and
different values of $\rho$ can be independently chosen by individual
devices or even individual terminals, all based only on local
information, such as the primal and dual residuals of a given device and
its neighboring nets.

\section{Conclusion}
We have presented a fully decentralized method for dynamic network
energy management based on message passing between devices. Prox-average
message passing is simple and highly extensible, relying solely on peer
to peer communication between devices that exchange energy.
When the resulting network optimization problem is convex, 
prox-average message passing converges
to the optimal value and gives optimal locational marginal prices. We
have presented a parallel implementation that shows the time per
iteration and the number of iterations needed for convergence of
prox-average message passing are essentially independent of the size of
the network. As a result, prox-average message passing can scale to
extremely large networks with almost no increase in solve time.

\section*{Acknowledgments}
The authors thank Yang Wang for extensive discussions on the problem
formulation as well as ADMM methods; Yang Wang and Brendan O'Donoghue for
help with the $\rho$ update method; and Ram Rajagopal and Trudie Wang for
helpful comments.

This research was supported in part
by Precourt 1140458-1-WPIAE, %TomKat
by AFOSR grant FA9550-09-1-0704, %MURI with Shanhui Fan
by AFOSR grant FA9550-09-0130, %Hearn AFOSR
and by NASA grant NNX07AEIIA. %NASA

\newpage
\bibliography{decen_dyn_opt}

\end{document}